\documentclass{article}
\usepackage{fullpage}
\usepackage{latexsym}
\usepackage{epsfig}
\usepackage{rotating}
\usepackage{amssymb}

\newcommand{\epsf}[1]{\epsfbox{#1}}

\newcommand{\GL}{\mathbb{L}}
\newcommand{\ns}{\mathbb{N}}
\newcommand{\zs}{\mathbb{Z}}

\newcommand{\qs}{\mathbb{Q}}
\newcommand{\rs}{\mathbb{R}}
\newcommand{\bx}{\bar x}
\newcommand{\by}{\bar y}
\newcommand{\bz}{\bar z}
\newcommand{\cs}{\mathbb{C}}
\newcommand{\cf}{\mathbb{F}}

\newcommand{\GE}{\mathbb{E}}

\newcommand{\Ref}[1]{(\ref{#1})}
\newcommand{\beq}{\begin{equation}}
\newcommand{\eeq}{\end{equation}}
\newcommand{\gf}{generating function}
\newcommand{\gfs}{generating functions}
\newcommand{\al}{\alpha}

\def\cqfd{\par\nopagebreak\rightline{\vrule height 3pt width 5pt depth 2pt}
\medbreak}

 \newtheorem{Theorem}{Theorem}
 \newtheorem{propo}[Theorem]{Proposition}
\newtheorem{coro}[Theorem]{Corollary}
\newtheorem{Lemma}[Theorem]{Lemma}

\title{\bf Walks on the slit plane  }
\author{
\parbox{8cm}{
\centerline {\sc Mireille Bousquet-M{\'e}lou\thanks{Both authors were partially
supported by the INRIA, via the cooperative research action
\textsf{Alcophys}.}}
\centerline{\small CNRS, LaBRI, Universit\'e Bordeaux 1}
\centerline{\small 351 cours de la Lib\'eration}
\centerline{\small 33405 Talence Cedex,  France}
\centerline{\small\tt bousquet@labri.u-bordeaux.fr }
}
\parbox{65mm}{
\centerline{\sc Gilles Schaeffer}
\centerline{\small CNRS, LORIA,
Campus Scientifique}
\centerline{\small 615 rue du Jardin Botanique -- B.P. 101}
\centerline{\small 54602 Villers--l\`es--Nancy Cedex,
France}
\centerline{\small\tt Gilles.Schaeffer@loria.fr }
}
}

\begin{document}
\maketitle

\begin{abstract}

In the first part of this paper, we enumerate exactly walks
on the square lattice that start from the origin, but otherwise avoid
the half-line ${\cal H}=\{(k,0), k\le 0\}$. We call them {\em walks on
the slit plane\/}. We count them by their
length, and by the coordinates of their endpoint. The corresponding
three variable \gf \ is algebraic of degree $8$. Moreover, for any
point  $(i,j)$, the length
\gf \ for walks of this type ending at  $(i,j)$ is
also algebraic, of degree $2$ or $4$, and involves the famous
Catalan numbers. 

Our method is based on the 
solution of a functional equation,
established via a simple combinatorial argument. It actually works for
more general models, in which walks take their steps in a finite
subset of $\zs ^2$ satisfying two simple conditions. The corresponding
\gfs \ are always algebraic.

In the second part of the paper, we derive from our enumerative
results a number of probabilistic corollaries. For instance, we can
compute exactly the probability that an ordinary random walk starting from
$(i,j)$ hits for the first time the half-line $\cal H$ at position
$(k,0)$, for any triple $(i,j,k)$. This generalizes a question raised
by R. Kenyon, which was the starting point of this paper. 

Taking uniformly at random all $n$-step walks on the slit plane, we
also compute the probability that they visit a given point $(k,0)$,
and the average number of visits to this point. In other words, we
quantify the transience of the walks. Finally, we derive an explicit
limit law for the coordinates of their endpoint.

\end{abstract}

\section{Introduction}
In January 1999, Rick Kenyon posted on the ``domino'' mailing-list the
following e-mail:
\begin{quotation}
``Take a simple random walk on $\zs^2$ starting on the $y$-axis
at $(0,1)$, and stopping when you hit the nonpositive
$x$-axis. 
Then the probability that you end at the origin is $1/2$.\\
Since this result was obtained from a long calculation
involving irrational numbers, I wonder if there is
an easy proof? By way of comparison,
if you start at $(1,0)$ the probability of stopping
at the origin in $2-\sqrt{2}$.''
\end{quotation}
This mail led  Olivier Roques, a graduate student at La{\sc BRI}, to investigate
the {\em number\/} of such walks of given length: he soon conjectured
that exactly $4^nC_n$ walks of length $2n+1$ go from $(0,1)$ to
$(0,0)$ without hitting the nonpositive $x$-axis before they reach
their endpoint, where $C_n={2n\choose n}/(n+1)$ is the $n$th Catalan
number.  This seems to confirm the statement that Catalan numbers  ``are
probably the most
frequently occurring combinatorial numbers after the binomial
coefficients''~\cite{sloane}. Similarly, O.~Roques conjectured that, if the starting point is
chosen to be $(1,0)$, then the number of walks is even more
remarkable, being $C_{2n+1}$.
These conjectures directly imply Rick Kenyon's results.

\begin{figure}[ht]
\begin{center}
\epsfig{file=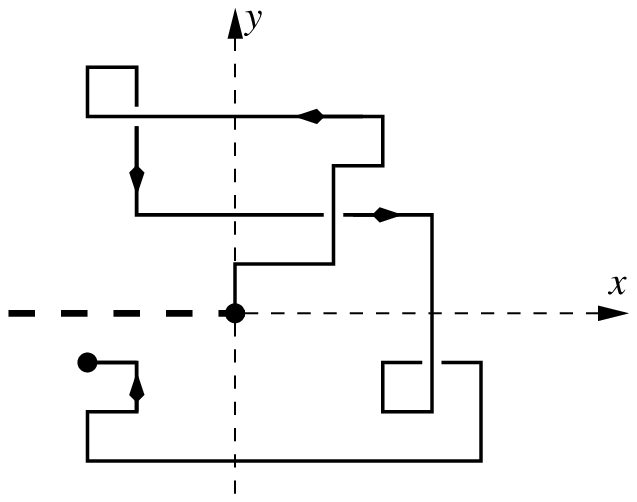}
\end{center}
\caption{A walk on the slit plane.}
\label{chemin}
\end{figure}

In this paper, we prove O. Roques's conjectures  as a special case
of a more complete result. More precisely, having
reverted the direction of the walks, we study the number
$a_{i,j}(n)$  of walks of length $n$ 
that start from $(0,0)$, end at $(i,j)$, and never hit the
horizontal half-axis ${\cal H}=\{(k,0): k\le 0\}$ once they have left their
starting point: we call them {\em walks on the slit plane\/}
(Fig. \ref{chemin}).
We give a closed form expression for the
{\em complete \gf \ }
$$
S(x,y;t) = \displaystyle 
\sum_{n\ge 0} \sum_{i \in \zs} \sum_{j\in\zs} 
a_{i,j}(n)x^iy^jt^n, $$
which turns out to be algebraic of degree $8$ (Theorem
\ref{dim2}). The series $S(1,1;t)$, which counts walks by their
length, has already been considered in the literature, and some
asymptotic estimates for its coefficients have been
obtained (see Lawler~\cite[Chap.~2]{lawler}). However, to our knowledge, it was
never realized that this walk model was exactly solvable. 
Note that a  refinement/variation
of an argument of Lawler is actually the starting point of another
possible derivation of $S(x,y;t)$ (see~\cite{prep}). 

Although O. Roques's conjectures beg for bijective
proofs, our approach is far from bijective. It is based on a functional
equation for $S(x,y;t)$,  which is extremely simple to establish: its
combinatorial counterpart is the deletion of the last step of the
walk. The originality of our method lies in the solution of this
equation, which might, at first, 
look a bit miraculous, and which
we believe to be of independent interest. Indeed, after staring at the
``miracle'' for a while, we realized that the same approach would
work just as well for a whole range of similar walk models where
the walks take their steps in a 
finite subset of $\zs ^2$ satisfying two simple conditions 
(the half-line $\cal H$ being still forbidden). These
conditions are satisfied, for example, by walks with diagonal steps
(Fig.~\ref{turndiag}). 
In all cases, the complete \gf \ for walks on the slit plane is
algebraic and can be written explicitly. This is  the central enumerative
result of our paper (Section~\ref{complete}).

Starting from this result, we then embark on the enumeration of walks on the
slit plane that end at a specified point $(i,j)$. Again, for all the
models that fit in our framework, these walks have an algebraic \gf \ 
(Section~\ref{prescribed}).  
In particular, for walks on the square lattice ending at $(0,1)$ or
$(1,0)$, we prove the Catalan results conjectured by O. Roques. 
Remarkably, walks with {\em diagonal\/} steps
also involve Catalan numbers. For instance, the number of walks ending
at $(1,1)$  is again $C_{2n+1}$. 
Note that ordinary walks on the (unslit) square lattice ending at
$(i,j)$ do {\em not\/} have an algebraic \gf . Nor do walks that
completely avoid the line $y=0$, and end at $(i,j)$.

Section~\ref{section-starting} is still devoted to exact enumeration,
but does not really contain original techniques: mostly, we need to
know, for the probabilistic applications we have in mind, the \gf \
for walks on the 
slit plane that start from the point $(k,0)$, and this can be derived
from the results previously obtained.

\medskip
So much for counting. The rest of the paper uses our enumerative
results to solve probabilistic questions, and focusses on the square lattice.
In Section~\ref{kenyon-like},  we generalize R. Kenyon's
problem by computing the probability that a random walk, starting from 
$(i,j)$, meets for the first time the half-line $\cal H$ at
abscissa $-k$. In other words, we compute exactly the {\em hitting
distribution\/} of $\cal H$, starting from $(i,j)$. Our results refine
asymptotics results from~\cite{kesten} and~\cite{lawler}.

In Section~\ref{section-transient}, walks on length $n$ on the slit
plane are taken uniformly at random, and we compute the probability $p_k$
that they visit the point $(k,0)$, in the limit $n \rightarrow
\infty$. As can be expected, the forbidden 
half-line creates a long-range repulsion, and the walks are {\em
transient\/}. That is, $p_k$ is strictly less than $1$. We also
compute the average number of visits to the point $(k,0)$, which is finite.

In Section~\ref{limit-laws},
walks of length $n$ are still taken uniformly at random, so that the position
$(X_n,Y_n)$ of their endpoint becomes a two-dimensional random
variable. We prove that $(X_n /\sqrt{n}, Y_n /\sqrt{n})$ converges in
distribution to an explicit law. 
In particular, the
average abscissa 
$\GE(X_n)$ grows like $c\sqrt{n}$, where
$c=\Gamma(3/4) /\Gamma(1/4)$: this confirms the fact that walks are
repelled away from the origin.

\bigskip
\noindent{\bf Notations.}  We shall use the following standard
conventions throughout the paper. Given a ring $\GL$ and 
$k$ indeterminates $x_1, \ldots , x_k$, we denote by

$\bullet$ $\GL[x_1, \ldots , x_k]$ the ring of polynomials in $x_1,
\ldots , x_k$ with coefficients in $\GL$,

$\bullet$ $\GL[[x_1, \ldots , x_k]]$ the ring of formal power series in $x_1,
\ldots , x_k$ with coefficients in $\GL$, that is, formal sums
$$\sum_{n_1\ge 0, \ldots , n_k \ge 0} a_{n_1, \ldots ,
n_k}x_1^{n_1}\cdots x_k^{n_k},$$ \\
and if $\GL$ is a field, we denote by

$\bullet$ $\GL(x_1, \ldots , x_k)$ the field of rational functions in $x_1,
\ldots , x_k$ with coefficients in $\GL$.\\
A {\em Laurent polynomial\/} in the $x_i$ is a polynomial in the
$x_i$ and the $\bx _i=1/x_i$. A {\em Laurent series\/} in the $x_i$ is
a series of the form $\sum_{n_1\ge N_1, \ldots , n_k \ge N_k} a_{n_1, \ldots ,
n_k}x_1^{n_1}\cdots x_k^{n_k},$ 
where $N_1, \ldots , N_k \in \zs$.


%
\section{The complete \gf } 
\label{complete}
\subsection{The ordinary square lattice}
We consider  walks on the square lattice made of
four kinds of steps: north, east, south and west.  Let $n\ge 0$, and
$(i,j) \in \zs\times \zs$.
We denote by $a_{i,j}(n)$ the number of walks of length $n$ 
that start from $(0,0)$, end at $(i,j)$, and never return to the
horizontal half-axis ${\cal H}=\{(k,0): k\le 0\}$ once they have left their
starting point: we call them {\em walks on the slit plane\/}.
Fig. \ref{chemin} shows such a walk, with $(i,j)=(-3,-1)$ and $n=46$.
We denote by $a(n)$ the total number of walks of length $n$ 
on the slit plane, 
regardless of their endpoint. 

Let $S(x,y;t)$ be the {\em complete generating function\/} for walks
on the slit 
plane, counted by their length and the position of their endpoint:
\begin{eqnarray}
S(x,y;t) & =&\displaystyle 
 \sum_{n\ge 0}  \sum_{i \in \zs} \sum_{j\in\zs}
a_{i,j}(n)x^iy^jt^n, \label{complete-definition} \\ 
&=& 1+t(x+y+\bar y)
+t^2(x^2+2xy+2x\bar y+\bar xy+\bx \by +y^2+\by ^2) \nonumber\\ 
&+& (5x +x^3 +4y + 4 \by +y^3 +\by ^3 +3xy^2 +3x\by^2 +3x^2y +3x^2\by
+\bx ^2 y+\bx ^2 \by
+2\bx y^2
+2\bx \by ^2)t^3+O(t^4) \nonumber
\end{eqnarray}
with the notations $\bar x=x^{-1}, \bar y =y^{-1}.$
This is a formal power series in $t$ with coefficients in $\rs[x, \bx]$.
We shall prove that this series is algebraic of degree $8$ over the
field of rational functions in $x,y$ and $t$, and actually give an
explicit expression for it.

\begin{Theorem}\label{dim2}
The complete generating function 
for walks on the slit plane is
$$
S(x,y;t) =\frac{\left(1-2t(1+\bar x)+\sqrt{1-4t}\right)^{1/2}
                   \left(1+2t(1-\bar x)+\sqrt{1+4t}\right)^{1/2}}
             {2(1-t(x+\bar x+y+\bar y))}.
$$
This series is algebraic of degree $8$. When $x=y=1$, it specializes
to
$$S(1,1;t)= \sum_{n \ge 0} a(n)t^n = 
\frac{(1+\sqrt{1+4t})^{1/2}(1+\sqrt{1-4t})^{1/2}}{2(1-4t)^{3/4}},$$
so that the asymptotic growth 
of the number of $n$-step walks on the slit plane is
$$a(n)\sim \frac{\sqrt{1+\sqrt{2}}}{2 \Gamma(3/4)} \ 4^n n^{-1/4}.$$
\end{Theorem}
In other words, the probability that a random walk  on the square
lattice, starting from $(0,0)$, has not  met the half-line $\cal H$
after $n$ steps is asymptotic to $cn^{-1/4}$ with $c= {\sqrt{1+\sqrt{2}}}/{2
/\Gamma(3/4)}$. The decay in $n^{-1/4}$ was known \cite[Eq.~(2.35)]{lawler},
but the detailed asymptotic behaviour of this probability
seems to be new.

We delay the proof of this theorem to state another result of the same
type. 

\subsection{The diagonal square lattice}
\begin{figure}[ht]
\begin{center}
\epsfig{file=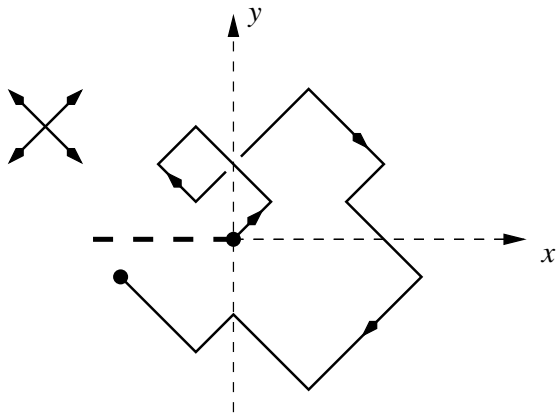}
\end{center}
\caption{A walk  on the slit plane with  diagonal steps.}
\label{turndiag}
\end{figure}
Let us consider the following variation on the previous model:
the forbidden half-line $\cal H$ is unchanged, but walks  now  consist of
\emph{diagonal} steps $(1,1),(1,-1),(-1,1)$ and $(-1,-1)$
(Fig.~\ref{turndiag}).  Let   $  a_{i,j}(n)$ be the
number of walks of length $n$ with diagonal steps that
start from $(0,0)$, end at $(i,j)$ (with $i+j$ 
even), and never return to $\cal H$. Let $S(x,y;t)$ denote the
corresponding complete \gf , defined as
in~\Ref{complete-definition}. We obtain for this series a result that 
is very similar to Theorem~\ref{dim2}, and
actually a bit simpler. 

\begin{Theorem}\label{dim2bis}
 The complete generating function 
for walks on the slit
plane  with diagonal steps is 
$$S(x,y;t)=
\frac{\left(1-8t^2(1+\bx ^2)+\sqrt{1-16t^2}\right)^{1/2}}
{\sqrt{2}(1-t(x+\bx)(y+\by))}.$$ 
This series is algebraic of degree $4$. When $x=y=1$, it specializes
to
$$S(1,1;t)= \sum_{n \ge 0} a(n)t^n = 
\frac{(1+4t)^{1/4}(1+\sqrt{1-16t^2})^{1/2}}{\sqrt{2}(1-4t)^{3/4}},$$
so that the asymptotic growth of the number of walks of length $n$ on
the diagonal slit plane is
$$a(n)\sim \frac{1}{2^{1/4} \Gamma(3/4)} \ 4^n n^{-1/4}.$$
\end{Theorem}
Theorems \ref{dim2} and \ref{dim2bis} look suspiciouly like each
other; we are going to show that they are two instances of a more
general result that applies to walks with rather general steps.

\subsection{A general result for walks on the slit plane}
\label{section-complete-general}
Let $A$ be a finite subset of $\zs ^2$. A walk with steps in $A$ is a
finite sequence $w = (w_0, w_1, \ldots , w_n)$  of vertices of $\zs
^2$ such that $w_{i}-w_{i-1} \in A$ for $1 \le i \le n$. The number of
steps, $n$, is the {\em length\/} of $w$. We say that $w$ avoids the
half-line ${\cal H}=\{(k,0), k\le 0\}$ if none of the vertices $w_1, \ldots ,
w_n$ belong to $\cal H$.
 The starting point of the walk, $w_0$, is  allowed to be on
$\cal H$. 
For $(i,j) \in \zs$ and $n \ge 0$, let
$a_{i,j}(n)$  denote the number of walks of length $n$, with steps in $A$,
that start from $(0,0)$, end at $(i,j)$, and avoid $\cal H$: we call them
{\em walks on the slit plane\/}. Let $S(x,y;t)$ be the associated
 {\em complete \gf }:  
$$
S(x,y;t) 
= \sum_{n\ge 0} \sum_{i \in \zs} \sum_{j\in\zs} 
a_{i,j}(n)x^iy^jt^n. $$
In what follows, we shall often omit the length variable $t$,
writing, for instance, $S(x,y)$ instead of $S(x,y;t)$.

\begin{figure}[ht]
\begin{center}
\epsfig{file=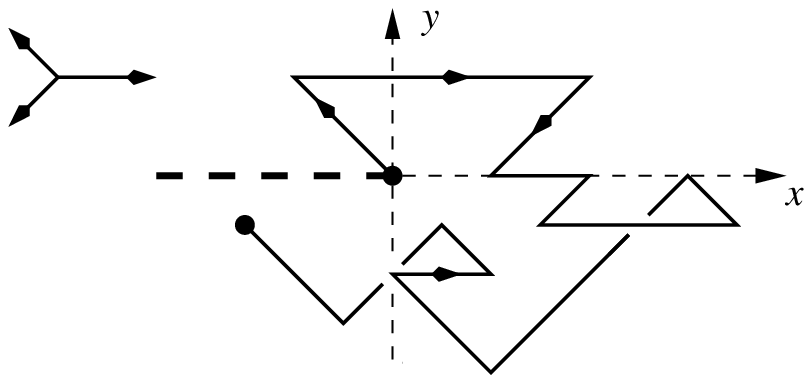}
\end{center}
\caption{A walk on the slit plane with steps in  $\{(-1,1),(-1,-1),(2,0)\}$.}
\label{triangulaire}
\end{figure}

We obtain a functional equation for the series $S(x,y)$ by saying
that a walk of length $n$  is
obtained by adding a step to another walk of length
$n-1$. However,  this procedure sometimes
produces a {\em bridge\/},  that is, a nonempty walk that 
starts  at
$(0,0)$, ends on the half-line $\cal 
H$, but otherwise avoids $\cal H$. Hence, denoting by $B(\bx)$ the \gf \ for
bridges, we have:
$$S(x,y)=1+tS(x,y) \left(\sum _{(i,j) \in A} x^i y^j \right)
-B(\bx),$$
that is,
\beq K(x,y) S(x,y)=1- B(\bx)
,\label{eqS-general}\eeq 
where 
\beq  K(x,y)=1-t\sum _{(i,j) \in A} x^i y^j \label{kernel-general}\eeq
is the {\em kernel\/} of  Eq.~\Ref{eqS-general}.
%
It turns out that  Eq.~\Ref{eqS-general}  can
be solved by elementary algebraic methods, provided
the set $A$ of allowed steps satisfies the two following conditions:
\begin{enumerate}
\item {\bf Symmetry}: the set of steps is symmetric with respect to the line
$y=0$; that is, if $(i,j) \in A$, then  $(i,-j) \in A$.
\item {\bf Small height variations}: for all $(i,j) \in A$, $|j|\le 1$.
\end{enumerate}
 From now on, we shall
restrict the study to sets $A$ satisfying these conditions.
Examples  include the ordinary
square lattice of Theorem~\ref{dim2}, the diagonal square lattice of
Theorem~\ref{dim2bis}, and the (oriented)
triangular lattice of Fig.~\ref{triangulaire}, corresponding to 
$A=\{(-1,1),(-1,-1),(2,0)\}$. The latter model is
equivalent, by  the transformation $(i,j) \rightarrow ((i+j)/2,j)$, to
the case 
$A=\{(0,1),(1,0),(-1,-1)\} $ recently studied by Ira Gessel via a
completely different approach~\cite{gessel}.

We define two Laurent polynomials in $x$:
\beq A_0(x)= \sum_{(i,0) \in A} x^i \ \ \ \hbox{and } \ \ \ 
A_1(x)=\sum_{(i,1) \in A} x^i.\label{A01}\eeq
We could define similarly a polynomial $A_{-1}(x)$, but the symmetry
condition makes it identical to $A_1(x)$. The kernel, given
by~\Ref{kernel-general}, can be rewritten 
\beq K(x,y)=1-tA_0(x)-t(y+\by) A_1(x).\label{kernel-A}\eeq
One of the main tools in the solution of this rather general model is
the following lemma, which defines a canonical way of factoring
polynomials in $t$ with coefficients in $\rs [x, \bx]$.
\begin{Lemma}[The factorization lemma]
\label{factorisation-lemma}
Let $\delta(x;t)$ be a polynomial in $t$  with coefficients in $\rs
[x, \bx]$, and assume  $\delta(x;0)=1$. There exists a unique triple
$(D(t),\Delta(x;t), \bar \Delta(\bx;t))\equiv(D,\Delta(x), \bar
\Delta(\bx))$ of 
formal power series in $t$ satisfying the following conditions: 
\begin{itemize}
\item $\delta(x) = D \Delta(x)\bar \Delta(\bx),$
\item the coefficients of $D$ belong to $\rs$,
\item the coefficients of $\Delta(x)$  belong to $\rs[x]$,
\item the coefficients of $\bar \Delta(\bx)$  belong to $\rs[\bx]$,
\item  $D(0)= \Delta(0;t)=\bar \Delta(0;t)=\Delta(x;0)=
 \bar \Delta(\bx;0)=1$.
\end{itemize}
Moreover,  these three series are algebraic, and  $\Delta(x)$ (resp.
$\bar \Delta(\bx)$) is a polynomial in $x$ (resp. $\bx$). 
\end{Lemma}
This lemma will be proved in Section~\ref{proofs-section}.
We can now state the main theorem of this section, which gives a closed
form expression for the complete \gf \ of walks on the slit plane.
\begin{Theorem}\label{dim2-general}
Let $A$ be a finite subset of $\zs ^2$ satisfying the conditions of
symmetry and small height variations.
Let $\delta(x)$ be the following polynomial
in $t,x$ and $\bx$:
\beq \delta(x) = \left( 1-t A_0(x) -2t A_1(x)\right)
\left( 1-t A_0(x) +2t A_1(x)\right).\label{delta-general}\eeq
Let $(D,\Delta(x), \bar \Delta(\bx))$ be the factorisation of
$\delta(x)$, defined as in the lemma above. Then the complete \gf \ for
walks on the slit plane with steps in $A$ is
$$S(x,y;t)= \frac{\displaystyle \sqrt{D \bar \Delta(\bx)}}{K(x,y)}$$
where the kernel $K(x,y)$ is given by~{\em \Ref{kernel-A}}.
The \gf \ for walks on the slit 
plane ending on the line $y=0$ is
$$S_0(x;t) := \sum_{n\ge 0} \sum_{i \in \zs} a_{i,0}(n)x^it^n=
\frac{1}{\sqrt {\Delta(x)}}.$$ 
These series  are algebraic.
\end{Theorem}
As a by-product of  this theorem, we are able to enumerate
bridges.
We shall use this result in
Section~\ref{section-starting}. 
\begin{propo}
\label{bridges0}
The \gf \ for bridges, counted by their length ($t$) and  abscissa of
their endpoint ($x$), is:
$$B(\bx;t)=1-\sqrt{D \bar \Delta(\bx)}.$$
\end{propo}
In probabilistic terms, this proposition gives the joint distribution
of the time $T$ and 
position $X$ where a random walk with steps in $A$, starting from
$(0,0)$, returns for the first time to $\cal H$ (we assume that the
walk takes each step of $A$ with uniform probability $1/a$, with
$a=|A|$):
$$\GE\left( t^T x ^X\right):= \sum_{n ,k\ge 0} \Pr(T=n, X=-k) t^n \bx ^k= 
1-\sqrt{D(t/a) \bar \Delta(\bx;t/a)}.$$

\subsection{Proofs}
\label{proofs-section}
\subsubsection{Another functional equation}
\label{establish}

Instead of working directly with
Eq.~\Ref{eqS-general}, it is convenient to introduce the  series
$$T(x,y;t)\equiv T(x,y)= \sum_{i, j, n } a_{i,j}(n)x^iy^{|j|} t^n,$$
which, due to the symmetry
condition, contains as much information
as $S(x,y)$, and has only nonnegative exponents in $y$.
Again, one writes easily a functional
equation for the series $T(x,y)$:
\beq K(x,y) T(x,y)=1- B(\bx)
+t S_0(x) A_1(x)(y-\by),\label{eqT-general}\eeq 
the kernel $K(x,y)$ being the same as in Eq.~\Ref{eqS-general}.
The role of the last term of this equation is to correct the power of
$y$ for walks whose last step starts on the line $y=0$ and ends at 
 ordinate $-1$.

\subsubsection{The kernel method}
\label{kernel}
The next ingredient of the proof is the so-called kernel method. The
principle is to cancel the kernel $K(x,y)$ by an appropriate choice of
$y$, so as to
obtain certain relations between the series that occur on the right-hand
side of the equation. See \cite{FaIa,knuth} for early uses of this
method, and \cite{6gus,BoPe98,Malyshev} for more recent developments.

 The kernel $K(x,y)$ of our functional equation is given by~\Ref{kernel-A}.
  As a Laurent polynomial in $y$, it  has two
roots. One of them is a formal power series in $t$:
$$Y(x)= \frac{1-tA_0(x)-\sqrt{\delta(x)}}{2tA_1(x)}
=A_1(x)t+A_0(x)A_1(x)t^2+O(t^3),$$
where $\delta(x)$ is the polynomial defined in Theorem~\ref{dim2-general}.
The coefficient  of $t^n$ in this series lies {\em a priori\/} in $\qs
(x)$, but the equation
$$Y(x)=tA_0(x)Y(x)+tA_1(x)\left( 1+ Y(x)^2\right)$$
implies that it is actually a Laurent polynomial in $x$.
By the symmetry assumption,  $K(x,y)= K(x,\by)$, and the other root of
the kernel is $1/Y(x)$.

Observe that $T(x,Y(x))$ is a well-defined series belonging to 
 $\zs[x,\ \bar x][[t]]$: this comes from the fact that $T(x,y)$ has
 only nonnegative exponents in the variable $y$.
Let us replace $y$ by $Y(x)$ in Eq.~\Ref{eqT-general}: the kernel
vanishes, leaving:
\begin{eqnarray}
1- B(\bx)
& =& t S_0(x) A_1(x)\left( \frac{1}{Y(x)}
-Y(x)\right)\nonumber \\
&=&S_0(x)\sqrt{\delta(x)}.\label{eq-delta}
\end{eqnarray}
Here, the facts that $y$ only appears 
as  $(y-\by)$ in the right-hand
side of~\Ref{eqT-general}, and that the kernel $K(x,y)$ is 
symmetric in $y$ and $\by$, play a crucial role. 
These two properties 
parallel (and actually inspired) the two properties imposed on $A$.

\subsubsection{Separation of the positive and negative parts in $x$}
\label{positive}
 Let us assume, for the moment, the existence of the three series 
 $D,\Delta(x)$ and $ \bar \Delta(\bx)$. 
We replace
 $\delta(x)$ by  $D\Delta(x) \bar \Delta(\bx)$  in 
Eq.~\Ref{eq-delta}, and divide the resulting equation by $\sqrt{\bar
\Delta(\bx)}$. This gives:
\beq\frac{1}{\sqrt{\bar \Delta(\bx)}}\left(1- B(\bx)\right)
=S_0(x)\sqrt{ D\Delta(x)}.\label{separe-general}\eeq
The two sides of this equation  are 
formal power series in $t$ whose
coefficients are Laurent polynomials in $x$. But we observe that the
right-hand side  only contains nonnegative exponents of $x$, while the
left-hand side  only contains nonpositive exponents of $x$. Moreover,
by definition of $S_0$ and $\Delta$, we have $\sqrt{\Delta(0)}
S_0(0)=1$. Let us now extract from Eq.~\Ref{separe-general} the terms with a
positive exponent of $x$. We obtain:
$$0= \sqrt{ D}\left(\sqrt{\Delta(x)} S_0(x)-1\right),$$
that is,
$$S_0(x)=\frac{1}{\sqrt{\Delta(x)}}.$$
Going back to Eq.~\Ref{separe-general}, this tells us that 
$$  1- B(\bx) ={\sqrt{D\bar\Delta(\bx)}}. $$
This is exactly Proposition~\ref{bridges0}. But the left-hand side of this equation is  the right-hand side of
Eq.~\Ref{eqS-general}, which defines the complete \gf \ $S(x,y)$. Hence
$$ S(x,y)= \frac{\sqrt{D\bar \Delta(\bx)}}{K(x,y)}. $$
This completes the proof of Theorem~\ref{dim2-general} (assuming
Lemma~\ref{factorisation-lemma} is true).

\subsubsection{The factorization lemma}
Let us finally  prove Lemma~\ref{factorisation-lemma}.
 Let the smallest
 exponent of $x$ occurring in $\delta(x)$ be  $-m$. Then $P(x;t)=x^m
\delta(x)$  is a
 polynomial in $x$ and $t$ such that $P(x;0)=x^m$. As a polynomial in
$x$, $P$ has degree, say, $d$, and hence admits $d$ roots,
denoted  $X_1, \ldots , X_d$, which belong to the algebraic closure of
$\rs (t)$. 
By Newton's theorem, there exists an integer 
$n\ge 1$ such that all these roots can be written as Laurent series
in the variable    $z=t^{1/n}$ (see~\cite[p.~89]{abhyankar}). 
Assume that  exactly $k$ of these roots, say
 $X_1, \ldots , X_k$, are finite at $z=0$.
The other $d-k$ roots %
contain terms of the form $z^{-i}$, with $i>0$.
The polynomial  $P(x;t)$  can be factored as
$$P(x;t)= x^m\delta(x) = D \prod_{i=1}^k (x- X_i)  \prod_{i=k+1}^d
\left(1-\frac{x}{  X_i}\right),$$
where $D$ is, as  the $X_i$, an algebraic function of $z$. For $i>k$,
the series $1/X_i$ equals $0$ at $z=0$. Hence the condition
$P(x;0)=x^m$ implies that $k=m$, that $D(0)=1$, and that the finite
roots  $X_1, \ldots , X_m$ equal $0$ when $z=0$. Let
$$\Delta(x)\equiv \Delta(x;z)=\prod_{i=m+1}^d \left(1-\frac{x}{
X_i}\right) \ \ \ \hbox{and}\ \ \ 
 \bar \Delta(\bx) \equiv \bar \Delta(\bx;z)=\prod_{i=1}^m (1-\bx X_i).$$ Then   the  series
 $D,\Delta(x)$ and $ \bar \Delta(\bx)$ satisfy all the required
conditions, but two: we still need to prove that they are actually
series in $t$ (and not only in  $z=t^{1/n}$) with real (rather than
complex) coefficients.
By  taking logarithms, we
obtain
$$ \log \delta(x) = \log D  + \log  \Delta(x)+ \log \bar \Delta(\bx).$$
As  $\Delta(0;z)=\Delta(x;0)=1$, the series $\log  \Delta(x)$
is  a multiple of 
$x$. Similarly,  the series $\log  \bar \Delta(\bx)$ is a multiple of
$\bx$.
But $\log \delta(x)$ can be written in a {\em unique\/}  way as
$$\log \delta= L_0(z) + xL_+(x;z) + \bx L_-(\bx;z)$$
where $L_0$ is a series in $z$ with real coefficients, independent of
$x$,  $L_+$ is a series in $z$ with 
coefficients in $\rs[x]$ and $L_-$ is a series in $z$ with
coefficients in $\rs[\bx]$. Moreover, $L_0, L_+$ and
$L_-$ are actually series in $t$, and not only in $z$.
The above equation forces $ D = \exp (L_0)$,
$\Delta(x)=\exp( xL_+(x))$ and $\bar \Delta(\bx)=\exp( \bx L_-(\bx))$,
and proves simultaneously that the series $D,\Delta(x)$ and $ \bar
\Delta(\bx)$ are actually series in $t$ with real coefficients, and
that they are the unique 
triple satisfying the required conditions.
This completes the proof of Lemma~\ref{factorisation-lemma}.

\subsubsection{Proof of Theorem~\ref{dim2}}
\label{premier}
For the ordinary square lattice, $A_0(x)=x+\bx$ and $A_1(x)=1$. The
polynomial $ {\delta(x)}=\left( 1-t(x+\bx +2)\right )\left( 1-t(x+\bx
-2)\right )$ has four roots $X_i$, $1 \le i \le 4$, which are
quadratic functions of $t$ and can be computed explicitly. Let $C(t)$
denote the \gf \ for Catalan numbers:
$$C(t)=\frac{1-\sqrt{1-4t}}{2t}=\sum_{n \ge 0}\frac{1}{n+1}{2n \choose
n} t^n.$$
Then $X_1= C(t)-1$ and $X_2=1-C(-t)$ are the two roots that are finite
at $t=0$, and by symmetry of $\delta$ in $x$ and $\bx$, the two other
roots are 
\beq X_3= (C(t)-1)^{-1} = \frac{1-2t+\sqrt{1-4t}}{2t} 
\quad \hbox{and  } \quad  X_4=(1-C(-t))^{-1}=
\frac{1+2t+\sqrt{1+4t}}{2t}.\label{X34}\eeq
By the previous subsection, the
canonical factorization of $\delta$ is such that
\begin{eqnarray}
  \Delta(x)=\bar\Delta(x)&=&(1-x X_3^{-1})(1-x X_4^{-1})
\label{Delta-formel} \\ 
&=&\left( 1-x (C(t)-1)\right) \left(
1+x(C(-t)-1) \right).
\label{Delta-def}
\end{eqnarray}
Taking the coefficient of $x^2$ in the relation $\delta(x)=D\Delta(x)
\Delta(\bx)$ yields
\begin{eqnarray}
 D&=& t^2 X_3 X_4 \label{D-formel}\\
&=&(C(t) C(-t))^{-2},\label{delta-Delta} 
\end{eqnarray}
as $(C(t)-1) = tC(t)^2$.
Now by Theorem~\ref{dim2-general}, Eqs.~\Ref{Delta-formel} and~\Ref{D-formel},
$$S(x,y;t) = \frac{\sqrt{t^2  X_3 X_4(1-\bx X_3^{-1})(1-\bx
X_4^{-1})}}
{K(x,y)} = \frac{\sqrt{(t X_3-t\bx)(t X_4-t\bx)}}
{1-t(x+\bx+y+\by)}.$$  
Theorem~\ref{dim2} then follows from~\Ref{X34}. We 
 obtain the asymptotic
behaviour of $a(n)$ by examining the singularities of
$S(1,1;t)$. (We refer to~\cite{flajolet-odlyzko} for a description of 
techniques that lead from the position and nature of the singularities
of a series to asymptotic estimates of its coefficients. These
techniques are applied in detail in Section~\ref{limit-laws}.)

\subsubsection{Proof of Theorem~\ref{dim2bis}}
\label{deuxieme}
For the diagonal square lattice,  $A_0(x)=0$ and $A_1=x+\bx$. Again, the
polynomial $ {\delta(x)}= 1-4t^2(x+\bx)^2$   has four roots $X_i$, $1 \le i \le 4$, which are quadratic
functions of $t$ and can be expressed explicitly in terms of the
Catalan \gf .  More precisely, $X_{1,2} =\pm 2t C(4t^2)$ and
\beq X_{3,4}=\pm X_1^{-1} = \pm\frac{1+\sqrt{1-16t^2}}{4t} 
.\label{X34-bis} \eeq
We then follow the same steps as above. The formal
expression of $\Delta(x)$  in
terms of the roots $X_i$ remains unchanged,  but the actual
value is of course different:
\begin{eqnarray} 
\Delta(x) =\bar\Delta(x) &=& (1-x X_3^{-1})(1-x X_4^{-1})\label{Delta-bis-formel}\\
&=&{1-4t^2x^2C(4t^2)^2}.
\label{Delta-bis}
\end{eqnarray}
 Taking the coefficient of $x^2$ in the relation $\delta(x)=D\Delta(x)
\Delta(\bx)$ yields now
\begin{eqnarray}
 D&=& -4t^2 X_3 X_4 \label{D-formel-2}\\
&=& C(4t^2)^{-2}.\label{delta-Delta-2} 
\end{eqnarray}
Theorem~\ref{dim2bis} then follows from  Theorem~\ref{dim2-general},
Eqs.~\Ref{Delta-bis-formel}, \Ref{D-formel-2}
and finally~\Ref{X34-bis}. 

\subsection{Refining the enumeration}
\label{section-refined}
In what we have done so far, each step of a walk contributes for a
weight $t$ in the \gf . We can refine the enumeration by giving
different weights to different steps, and the method works just as
well. Let us, for instance,  take into account the
number of vertical steps in the first model (ordinary square
lattice). 
Let us denote by  $S(x,y; t,v)$ the refined generating function, in which
 the variable $v$ keeps track of the vertical steps.
The polynomial $\delta$ is now
$$ \delta(x) = (1-t(x+\bx+2v))(1-t(x+\bx-2v)).$$
The computations then follow exactly the same lines as in
Section~\ref{premier}. The roots of $\delta$ that diverge when $t=0$ are
$$ X_3=\frac{1-2tv+\sqrt{\delta_1}}{2t}, 
\quad X_4= \frac{1+2tv+\sqrt{\delta_2}}{2t},
$$
with
\beq \delta_1=(1-2t(1+v))(1+2t(1-v))
\ \ \ \hbox{ and }\ \ \  \delta_2=
(1-2t(1-v))(1+2t(1+v)).\label{delta12}\eeq
We still have 
\beq \Delta(x)=\bar\Delta(x)=(1-x X_3^{-1})(1-x X_4^{-1}) \label{Delta-v}\eeq
and $D=t^2X_3X_4$.
%
We  thus obtain the following refinement of Theorem~\ref{dim2}.
\begin{Theorem}
\label{dim2-refined}
On the ordinary square lattice, the refined \gf \ $S(x,y;t,v)$ for
walks on the slit plane is  given by 
$$
S(x,y;t,v) =\frac{\left(1-2t(v+\bar x)+\sqrt{\delta_1}\right)^{1/2}
                   \left(1+2t(v-\bar x)+\sqrt{\delta_2}\right)^{1/2}}
             {2\left(1-t(x+\bar x+yv+\bar y v)\right)},
$$
where $\delta_1$ and  $\delta_2$ are given by~{\em \Ref{delta12}}.
\end{Theorem}

\section{Walks ending at a prescribed position}
\label{prescribed}
We still consider a model of walks on the slit plane where the set $A$
of allowed steps satisfies the conditions of symmetry and small height
variations that led to the expression of the complete \gf \ $S(x,y;t)$
given in Theorem~\ref{dim2-general}. 
We would like to enumerate walks on the slit plane ending at a
prescribed point $(i,j)$. 
Let 
$$S_{i,j}(t) \equiv S_{i,j} = \sum_{n \ge 0} a_{i,j}(n)t^n$$
be the corresponding \gf . This series is obtained by extracting the
coefficient of $x^iy^j$ from $S(x,y;t)$. As $i$ and $j$ belong to
$\zs$, rather than $\ns$, this is not an obvious task. In particular,
the algebraicity of $S(x,y;t)$ does {\em not\/}
automatically imply the  algebraicity of $S_{i,j}(t)$.
This is clearly shown by the enumeration of walks starting
from $(0,0)$ in the ordinary, unslit square lattice. The complete \gf \ is 
 $1/[1-t(x+\bx+y+\by)]$. It is rational, hence algebraic. However, for
$i$ and $j$ in $\zs$, the coefficient of $x^iy^j$ in this series is
$$\sum_{n \ge 0} {{n} \choose \frac{n+i+j}{2}}{ {n
} \choose \frac{n-i+j}{2}} t^n,$$
where the sum is restricted to integers $n$ of the same parity as
$i+j$. This series is transcendental: the coefficient of $t^n$ grows
like $ 4^n/n$, up to a multiplicative constant, revealing a logarithmic
singularity in the \gf \ that implies its transcendence (see
\cite{Flajolet87} for a discussion on the possible singularities of an
algebraic series).

In contrast, we shall prove that for all models covered by
Theorem~\ref{dim2-general},  the series
$S_{i,j}(t)$
is algebraic for all $i$ and $j$. We shall  describe, in terms of
the set $A$ of steps,  an algebraic
extension of $\qs(t)$ that contains these series. Our approach is
constructive, and gives a procedure that computes $S_{i,j}$
explicitly, for a given pair $(i,j)$. 
We have implemented this procedure as a \texttt{Maple} program.
Using these programs, we were able,  for instance, to compute
 $S_{i,j}$ for $|i|+|j|\leq 10$ for the above mentioned
square, diagonal and triangular lattices.

\subsection{The ordinary and diagonal square lattices}
For the two models that we studied in detail in
Section~\ref{complete},  the series  $S_{i,j}$ admit a rational
expression in terms of 
the following power series in $t$: 
\beq 
u=\frac{\sqrt{1+4t}-1}
{\sqrt{1-4t}+1}= \sum_{n\ge 0} (2.4^nC_n
-C_{2n+1})t^{2n+1}, \label{u-def}\eeq 
where $C_n=\frac{1}{n+1}{2n\choose n}$ is the $n$th Catalan number. 
Note that $u$ is quartic over $\qs(t)$:
$$t=\frac{u(1-u^2)}{(1+u^2)^2}.$$
This equation allows us to write any rational function of $t$ as a
rational function of $u$, and implies that
$\qs(u)=\qs(t,\sqrt{1-4t},\sqrt{1+4t})$. 
 We shall also use the following identities, relating 
the Catalan generating function $C(t)$ to the series $u$:
\beq C(t)=\frac{1+u^2}{1-u}, \ \ \ 
C(-t)=\frac{1+u^2}{1+u} , \ \ \
 tC(t)C(-t)=u,
\ \ \ \hbox{ and } \ \ \
C(4t^2)=\left( \frac{1+u^2}{1-u^2}\right)^2. \label{Catalan-u}\eeq

\begin{Theorem}[Ordinary square lattice: Prescribed endpoint]\label{dim0}
 For all $i$ and $j$,  the  \gf \ $S_{i,j}(t)=\sum_n
a_{i,j}(n)t^n$  for walks on the slit plane
ending at $(i,j)$ belongs to $\qs(u)=\qs(t,\sqrt{1-4t},
\sqrt{1+4t})$ and can be computed explicitly. It is either
 quadratic, or quartic over $\qs(t)$.
In particular,
$$S_{0,1}(t) =
\frac{u}{1-u^2}=
\frac{1-\sqrt{1-16t^2}}{8t}=\sum_{n\ge0} 4^nC_{n}t^{2n+1},$$ 
and
$$S_{1,0}(t) 
=\frac{u(1+u^2)}{1-u^2}=
\frac{2-\sqrt{1-4t}-\sqrt{1+4t}}{4t}= \sum_{n\ge 0} C_{2n+1}t^{2n+1},$$
as conjectured by O. Roques.
Some other  values are
$$S_{-1,1}(t)
= \frac{u^2}{1-u^2}=
\frac{\sqrt{1+4t}-\sqrt{1-4t}-4t}{8t}= \frac{1}{2} \sum_{n\ge 1}
C_{2n}t^{2n}, $$
$$S_{1,1}(t)
=\frac{u^2(2-u^2)}{(1-u^2)^2}=
\frac{1-24t^2+4t\sqrt{1+4t}-4t\sqrt{1-4t}-\sqrt{1-16t^2}}{32t^2}= \sum_{n\ge 1}
(4^{n-1}C_n+C_{2n}/2) t^{2n}.$$
\end{Theorem}
{\bf Conjecture.} In the early days of this study, it was observed by
B\'etr\'ema that the coefficients of $S_{-i,i}(t)$, for $i \ge 1$, seem
to factor nicely. The program we have written to compute explicitly
the series $S_{i,j}$ confirms this observation, and we conjecture
that, 
for $i \ge 1$ and  $n \ge i$,
$$a_{-i,i}(2n)= \frac{i}{2n} {{2i} \choose i}{{n+i}\choose {2i}}
\frac{{{4n}\choose{2n} }}{{{2n+2i}\choose {2i}}}.$$  The Lagrange
inversion formula (see~\cite[Sec.~1.2.4]{gj}), combined with a
binomial identity, implies that this conjecture is equivalent to
$$S_{-i,i}(t)= 
\frac{(-1)^i}{(1-u^2)^{2i-1}}
\sum_{k=i}^{2i-1}{2i-1\choose k}(-1)^ku^{2k}.
$$
We have verified this conjecture up to $i=10$.

\medskip
We now state a theorem, similar to Theorem~\ref{dim0}, for walks with
diagonal steps. Recall that, for this model, the \gf \ for walks
ending on the line $y=0$ is especially simple: 
$S_0(x;t)= (1-4t^2x^2 C(4t^2))^{-1/2}$. This allows us to write a
closed form expression for the series $S_{2i,0}$ and its
coefficients. By some aspects, the diagonal model looks sometimes simpler
than the ordinary one. By some aspects only: we shall see that
computing $S_{i,j}$, for $j \not = 0$, is significantly more difficult
for the diagonal model than for the ordinary one, even though the
final results look close.

\begin{Theorem}[Diagonal square lattice: Prescribed endpoint]
\label{dim0bis}  For all $i$ and $j$,  the  \gf \ $S_{i,j}(t)=\sum_n
a_{i,j}(n)t^n$  for walks on the slit plane
ending at $(i,j)$ belongs to $\qs(u)=\qs(t,\sqrt{1-4t},
\sqrt{1+4t})$ and can be computed explicitly. It is either
 quadratic, or quartic over $\qs(t)$.
In particular,
$$ S_{1,1}(t) 
=\frac{u(1+u^2)}{1-u^2}=
\frac{2-\sqrt{1-4t}-\sqrt{1+4t}}{4t}= \sum_{n\ge 0} C_{2n+1}t^{2n+1},$$
$$ S_{-1,1}(t)
=u=\frac{\sqrt{1+4t}-1}{\sqrt{1-4t}+1}=\sum_{n\ge 0}
(2.4^nC_n-C_{2n+1})t^{2n+1},$$
$$ S_{0,2}(t) 
=\frac{2u^2}{1-u^2}=\frac{\sqrt{1+4t}-\sqrt{1-4t}-4t}{4t}=
\sum_{n\ge 1} C_{2n} t^{2n},$$
and for $i \ge 1$,
$$
S_{2i,0}(t)
= {{2i} \choose {i}}  t^{2i}C(4t^2)^{2i}
=\sum_{n\ge i} \frac{i}{n}{ {2i}  \choose {i}}
{{2n} \choose {n-i}} 4^{n-i} t^{2n}.$$
\end{Theorem}
\noindent {\bf Remark.} Performing a counterclockwise rotation of 45 degrees on
Fig.~\ref{turndiag} shows that counting walks on the slit plane with
diagonal steps is equivalent to counting walks on the {\em ordinary\/}
square lattice that avoid the diagonal half-line $\{(k,k), k\le 0\}$. 
 Theorems~\ref{dim0} and~\ref{dim0bis} imply that,  among the walks
of length $2n+1$ that go from $(0,0)$ to $(1,0)$ on the ordinary
square lattice, exactly as many
avoid the horizontal half-line $\{(k,0), k\le 0\}$ as the diagonal
half-line $\{(k,k),k\le 0\}$ (and this number is $C_{2n+1}$). It would
be interesting to find a combinatorial explanation for this fact.

\medskip
Once again, Theorems~\ref{dim0} and \ref{dim0bis} follow from a more
general result that applies to walks with symmetry and small height
variations.

\subsection{Computing $S_{i,j}$ in the general case}
We start, naturally, from  the complete \gf \
$S(x,y;t)=\sum_{i,j \in \zs} S_{i,j}(t) x^i y^j$ 
derived in  Section~\ref{complete}. Our first step is fairly simple, and
consists in extracting  from
$S(x,y;t)$ the coefficient of $y^j$, 
denoted $S_j(x;t)$:
$$S_j(x;t)=\sum_{n  \ge 0} \sum_{i \in \zs} a_{i,j}(n) x^i t^n.$$
Of course, this series is the \gf \ for walks
ending at ordinate $j$. 
%
%
%
By symmetry, we can assume $j\ge 0$. The
polynomials $A_0(x), A_1(x)$ and $\delta(x;t)$ are defined as in
Section~\ref{section-complete-general} (see~\Ref{A01} and
\Ref{delta-general}). The canonical factorisation of 
$\delta(x)$ is still 
denoted $(D(t), \Delta(x;t), \bar\Delta(\bx,t))$.
\begin{Lemma}\label{fixed-ordinate-general}
For $j \ge 0$, the \gf \ for walks on the slit 
plane ending on the line $y=j$ is algebraic, and admits the following
expression: 
\beq 
S_{j}(x;t)=
\frac{1}{(2tA_1(x))^j}\left(
\frac{f^+_j(x;t)}{\sqrt{\Delta(x)}}-
{f^-_j(x;t)}{\sqrt{D\bar\Delta(\bar x)}}\right),
\label{almost-split} \eeq
where $f^+_j(x,t)$ and  $f^-_j(x,t)$ are the following polynomials in
$x, \bx$ and $t$:
\beq f_j^+(x;t)= \sum_{k=0}^{\lfloor j/2\rfloor}
{j \choose {2k}} \delta(x)^k (1-tA_0(x))^{j-2k},
\label{fj+}\eeq
\beq f_j^-(x;t)= \sum_{k=0}^{\lfloor (j-1)/2\rfloor}
{j \choose {2k+1}} \delta(x)^k (1-tA_0(x))^{j-2k-1}.\label{fj-}\eeq
\end{Lemma}
{\bf Proof.}
We start from the expression of $S(x,y;t)$ given in Theorem~\ref{dim2-general}.
Note that $1/K(x,y)$ is the \gf \ for unrestricted
walks on the plane. In Section~\ref{kernel}, we have seen that the two
roots of $K(x,y)$, seen as a polynomial in $y$, are the series $Y(x)$
and $1/Y(x)$, where 
$$Y(x)= \frac{1-tA_0(x)-\sqrt{\delta(x)}}{2tA_1(x)}.$$
  Let us convert the rational function $1/K(x,y)$ into partial
fractions of $y$:
\begin{eqnarray*}
\frac{1}{K(x,y)}&=& 
\frac{1}{\sqrt {\delta(x)}}
\left( \frac{Y(x)}{y- Y(x)}+ \frac{1}{1-yY(x)} \right) \\
&=& \frac{1}{\sqrt {\delta(x)}}
\left( \frac{\by Y(x)}{1- \by Y(x)}+ \frac{1}{1-yY(x)} \right).
\end{eqnarray*}
As $\delta(x) =1+O(t)$ and $Y(x)=O(t)$,  this identity splits
$1/K(x,y)$ as the sum of two power series in $t$ 
with coefficients in $\qs[x,\bx,y]$,  
one with only negative powers of $y$, and the other
with only nonnegative powers of $y$. 
Hence, by Theorem~\ref{dim2-general}, the coefficient of $y^j$ in
$S(x,y;t)$ is, for $j \ge 0$,  
\beq S_j(x)=  [y^j] S(x,y;t)=\sqrt{D\bar \Delta(\bx)} \ [y^j] \frac{1}{K(x,y)}
=\frac{\sqrt{D\bar \Delta(\bx)}}{\sqrt{\delta(x)}}Y(x)^j =
\frac{Y(x)^j}{\sqrt {\Delta(x)}}.\label{Sj-concis}\eeq
The binomial formula, applied to $Y(x)^j$, combined with the fact that
$\delta(x) = D \Delta(x) 
\bar \Delta(\bx)$, gives the announced expression of $S_j(x;t)$. This
expression is certainly much bigger than~\Ref{Sj-concis}, but also
more convenient for the purpose we have in mind. \cqfd

It is easy to see that, if a series $S(x;t)$ of $\qs [x][[t]]$ is
algebraic, then for $i \ge 0$, the coefficient of $x^i$ in this series
is a power series in $t$ that is also algebraic. 
(This can be proved by induction on $i$, by differentiating $i$ times
with respect to $x$ the polynomial equation satisfied by $S$ and then
setting $x=0$.)
However, this is not true, in general, if $S(x;t)$ belongs to $\qs
[x,\bx][[t]]$.
Therefore, in order to prove the algebraicity of the series
$S_{i,j}(t)$, our next step will be to separate in $S_j(x;t)$ the
terms with a negative power of $x$ from the other terms. In general,
given a series $S(x;t)$ belonging to $\qs [x,\bx][[t]]$, we define the
positive and negative parts of $S$ by:
$$S(x;t)=\sum_{n \ge 0} t^n \sum_{i \in \zs}a_i(n) x^i \Longrightarrow
 S^+(x;t)=\sum_{n \ge 0} t^n \sum_{i \ge 0}a_i(n) x^i \ \ \ \hbox{and }\
\ \  S^-(x;t)=\sum_{n \ge 0} t^n \sum_{i < 0}a_i(n) x^i.$$
 Let us examine Expression~\Ref{almost-split} of $S_j(x;t)$. 
As $\Delta(x)$ contains only positive powers of $x$, and $\bar
 \Delta(\bx)$ only negative powers, it is tempting to conclude that
 this expression 
essentially splits the series
 $S_j(x;t)$ into its  positive and negative parts. Indeed, as
 $f_j^+(x;t)$ and $f_j^-(x;t)$ are Laurent polynomials, it will always be easy
 to extract the positive  part of the portion of the
 expression that lies between the brackets: expanding  $1/\sqrt{\Delta(x)}$ in
$x$ allows us to remove from $f_j^+(x;t)/\sqrt{\Delta(x)}$ the
(finitely many) terms with a negative exponent of $x$, and similarly,
expanding  $\sqrt{\bar \Delta(\bx)}$ in 
$\bx$ allows us to remove from $f_j^-(x;t)\sqrt{\bar \Delta(\bx)}$ the
(finitely many) terms with a nonnegative exponent of $x$. 
To illustrate this procedure, let us consider the case $j=1$ of
 the ordinary square lattice:
$$S_1(x;t)= \frac{1}{2t} \left( \frac{1-t(x+\bx)}{\sqrt{\Delta(x)}}
-\sqrt{D \Delta(\bx)}\right),$$
where $\Delta(x)$ and $D$ are given by~\Ref{Delta-def} and~\Ref{delta-Delta}.
This readily gives
\begin{eqnarray}
S_1^+(x;t)= \sum_{i \ge 0, n\ge 0} a_{i,1}(n) x^i t^n &=&\frac{1}{2t} \left(
\frac{1-t(x+\bar x)}{\sqrt{\Delta(x)}}
+t\bx -\sqrt D\right),\label{S1+carre}\\
S_1^-(x;t)= \sum_{i < 0, n\ge 0} a_{i,1}(n) x^i t^n&=&
\frac{1}{2t}\left( \sqrt D -t\bx -\sqrt{D\Delta(\bx)}\right).\label{S1-carre}
\end{eqnarray}
As $A_1(x)=1$ for the ordinary square lattice, the extraction of the
positive part of~\Ref{almost-split} reduces for any $j$ to the
extraction of the positive part of the expression that lies between
the brackets.  We have just argued that this is an easy task, and we
can actually describe the form of $S_j^+(x)$. It is, essentially,
$f^+_j(x;t)/\sqrt{\Delta(x)}/(2t)^j$, plus some corrections: these
corrections are polynomials in $x$ and $\bx$, whose coefficients
involve $t$, $\sqrt D$, and the coefficients of $\Delta(x)$
seen  as a polynomial in $x$.

In general, however, problems will arise from nontrivial polynomials
$A_1(x)$.  For instance, $A_1(x)=x+\bx$ for the diagonal square
lattice model, and we need a special procedure to extract the positive
part of the series 
\beq S_1(x;t)= \frac{1}{2t(x + \bx )} \left(
\frac{1}{\sqrt{\Delta(x)}} -
\sqrt{D\Delta(\bx)}\right)\label{S1-diagonal}\eeq where $\Delta(x)$
and $D$ are now given by~\Ref{Delta-bis}
and~\Ref{delta-Delta-2}. 
The term $1/(2t(x+\bx)\sqrt{\Delta(x)})$, which is certainly a
tempting candidate for the positive part of $S_{1}$, is unfortunately 
a series in
$t$ with rational (rather than polynomial) coefficients in $x$. 

The main theorem of this section states that in general,  the positive part of
$S_j$  has the expected form  
(almost the left side of
Expression~\Ref{almost-split}). But the 
correction terms can be more subtle than for the ordinary
square lattice. They belong to 
an  algebraic extension of $\qs(t)$
which, as expected,  contains $\sqrt D$, the coefficients of
$\Delta(x)$ and $\bar 
\Delta (\bx)$ (seen as polynomials in $x$ and $\bx$ respectively), but
also the series $\sqrt{\Delta(\alpha _i)}$, where the numbers $\alpha _i$
are the roots of the polynomial $A_1(x)$, and the algebraic numbers
$\alpha _i$ themselves. The finite extension of $\qs(t)$  generated by
these algebraic numbers and
series, which only depends on the set $A$ of 
steps, is denoted $\qs_A$ below.
\begin{Theorem}[General model: Prescribed endpoint]
\label{General-endpoint}\label{dim1}
Let $A$ be a finite subset of $\zs ^2$ satisfying the conditions of
symmetry and small height variations. 
Let $j \ge 0$.
There exists a Laurent polynomial in $x$, with coefficients in
the field $\qs_A$ defined just above, denoted  $g_{j}(x;t)$,  such that
\begin{eqnarray*}
S_{j}^+(x;t)&=&\sum_{n\geq0} \sum_{i\ge 0}a_{i,j}(n)x^it^n
\;=\;
\frac{1}{(2tA_1(x))^j}\left(
\frac{f^+_j(x;t)}{\sqrt{\Delta(x)}}-g_{j}(x;t)\right),
\end{eqnarray*}
and
\begin{eqnarray*}
S_{j}^-(x;t)&=&\sum_{n\geq0} \sum_{i<0}a_{i,j}(n)x^it^n
\;=\;
\frac{1}{(2tA_1(x))^j}\left(g_{j}(x,t)-
{f^-_j(x;t)}{\sqrt{D\bar\Delta(\bar x)}}\right),
\end{eqnarray*}
where
$f_j^+(x;t)$ and $f_j^-(x;t)$ are the polynomials in $x, \bx$ and $t$
given by~{\em \Ref{fj+}} and~{\em \Ref{fj-}}.
Consequently, for any $(i,j)\in\mathbb{Z}^2$, the generating function 
$
S_{i,j}(t)\;=\;\sum_{n\geq0}a_{i,j}(n)t^n
$
for walks on the slit plane 
that end at position $(i,j)$ belongs to $\qs_A$. 
\end{Theorem}

\subsection{Proofs}
\subsubsection{The extraction procedure}
In order to prove  Theorem~\ref{General-endpoint}, we have to
learn how to extract the 
positive part of
certain series of the form $U(x;t)/A(x)$, where $A(x)$ is a
polynomial. The main tool will be a good old Taylor expansion.
\begin{Lemma} \label{taylor}
 Let $\alpha$ be a nonzero complex number, and $m$ a positive integer.
We define an operator $T_{\alpha, m}$  on
series of $\cs[x, \bx][[t]]$
by:
$$T_{\alpha, m} (U(x;t)) = \frac{1}{(x-\alpha)^m}\left( U(x;t) - \sum_{k=0}^{m-1}
U^{(k)}(\alpha;t) \frac{(x-\alpha)^k}{k!}\right),$$
where $U^{(k)}$ denotes the $k$th derivative of $U$ with respect to
$x$. 
If $U(x;t)$ belongs to $\cs[x,\bx][[t]]$ (resp.~$\cs [x][[t]]$, 
resp.~$\bx\cs [ \bx][[t]]$), then so does 
the series $T_{\alpha, m} (U(x;t))$.
\end{Lemma}
{\bf Proof.} The operator $T_{\alpha,m}$ acts 
coefficient-wise
on power series of $t$, so that we can concentrate on its effect
on $\cs[x,\bx]$. 
By linearity, it suffices to prove the statement when $U(x)= x^n$,
for $n \in \zs$. In this case, let
$$R(x):= U(x) - \sum_{k=0}^{m-1} U^{(k)}(\alpha)
\frac{(x-\alpha)^k}{k!}.$$ 
Taylor's formula implies that
$R^{\ell}(\alpha)=0$ for all $0\leq \ell < m$.  
When $R(x)$ is a
polynomial (\emph{i.e.} when $n\geq0$), this immediately implies its
divisibility by $(x-\alpha)^m$, and we conclude that
$T_{\alpha,m}(U(x))=R(x)/(x-\alpha)^m$ is a polynomial in $x$ as well.

Now if $U(x)=\bx ^n$ with $n >0$, then $P_1(x)=x^nR(x)$ is a
polynomial, with degree at 
most $n+m-1$. Using Leibnitz formula for successive derivatives, we
also have $P_1^{(\ell)}(\alpha)=0$ for $0\leq\ell< m$ and this implies
that $P_1(x)=(x-\alpha)^mP_2(x)$ with $P_2(x)$ of degree at most
$n-1$. Finally $T_{\alpha,m}(U(x))=\bar x^nP_2(x)$ belongs to $\bx\cs[\bx]$.
\cqfd
This lemma has a straightforward corollary which is a promising first step
towards the extraction of the positive part of $S_j(x;t)$.
\begin{coro}\label{diviseur-simple}
Let $\alpha$ be a 
complex number, and $m$ a positive integer. Let  $W(x;t) $ be  a
series of $\cs [x, \bx][[t]]$, and define $U(x;t)= (x-\alpha)^m
W(x;t)$. Clearly,  $U(x;t)$ also belongs to  $\cs [x,
\bx][[t]]$. Moreover,
$$W^+(x;t) = T_{\alpha, m} (U^+(x;t)).$$
(If $\alpha=0$, the operator $T_{\alpha,m}$ can be defined on
series, which, like $U^+(x;t)$, have their coefficients in $\cs[x]$.)
\end{coro}
{\bf Proof.} 
The result is obvious if $\al =0$. Otherwise, we derive from   the
fact that  $U(x;t)= (x-\alpha)^m 
W(x;t)$ that, for $0\le k <m$,
$U^{(k)}(\alpha;t) =0.$
Consequently,
\begin{eqnarray*}
W(x;t) &=& \frac{U(x;t)}{(x-\alpha)^m}
\;=\;  \frac{1}{(x-\alpha)^m}\left( U(x;t) - \sum_{k=0}^{m-1}
U^{(k)}(\alpha;t) \frac{(x-\alpha)^k}{k!}\right)\\
 &=& T_{\alpha, m} (U(x;t))\\
 &=& T_{\alpha, m} (U^-(x;t)) + T_{\alpha, m} (U^+(x;t)).
\end{eqnarray*}
The statement now follows  from Lemma~\ref{taylor}. 
\cqfd
There remains to iterate the above corollary to extract the positive
part of a series of the form $U(x;t)/A(x)$, where $A(x)$ is a
polynomial in $x$ with complex coefficients.
\begin{propo}\label{extract-positive}
Let $A(x)$ be a non-zero polynomial in $x$ with complex coefficients. Assume
$A(x)$ has exactly $d$ distinct roots, $\alpha_1,
\ldots , \alpha_d$, of multiplicities $m_1, \ldots , m_d$. Let
$W(x;t)$ be a series of $\cs [x, \bx][[t]]$, and define
$U(x;t)=A(x)W(x;t)$.  Clearly, $U(x;t)$ also belongs to $\cs [x,
\bx][[t]]$. There exist polynomials $P_{i,k}(x)$ with complex
coefficients such that
$$W^+(x;t)=  \frac{1}{A(x)}\left( U^+(x;t) -\sum_{i=1}^d \sum
_{k=0}^{m_i-1} P_{i,k}(x) U^{+(k)}(\alpha_i;t)\right) .$$
Moreover, if $A(x)$ and $W(x;t)$ have their coefficients in $\qs$ rather than
$\cs$, then the  polynomials
$P_{i,k}(x)$ have their coefficients in 
$\qs(\alpha_1, \ldots , \alpha_d)$.
\end{propo}
{\bf Proof.}
We proceed by induction on the number $d$ of distinct roots of
$A(x)$. If $d=0$, then $A(x)$ is simply a complex number, and the
result is obvious.
Now, assume the
result holds for $d-1$ distinct roots, with $d\ge 1$. Let $A(x)$ be a
polynomial with $d$ distinct roots  $\alpha_1, \ldots , \alpha_d$, of
multiplicity $m_1, \ldots , m_d$. 
For the sake of
simplicity, let us denote $\al \equiv \al_d$ and $m\equiv m_d$.
Let $V(x;t)=(x-\al)^m W(x;t)$. By Corollary~\ref{diviseur-simple},
\beq W^+(x;t) = \frac{1}{(x-\alpha)^m}\left( V^+(x;t) - \sum_{k=0}^{m-1}
V^{+(k)}(\alpha;t) \frac{(x-\alpha)^k}{k!}\right).\label{uneracine}\eeq
Now, writing $A(x)=(x-\al)^m B(x)$, where $B(x)$ is a polynomial with
only $d-1$ distinct roots, we have $U(x;t)=B(x)V(x;t)$, so that, by
the induction hypothesis,
\beq V^+(x;t)=  \frac{1}{B(x)}\left( U^+(x;t) -\sum_{i=1}^{d-1} \sum
_{k=0}^{m_i-1} Q_{i,k}(x) U^{+(k)}(\alpha_i;t)\right)
\label{d-1racines}\eeq
for some polynomials $ Q_{i,k}(x)$.
The expression of $W^+(x;t)$  follows, by combining~\Ref{uneracine}
and~\Ref{d-1racines}. 
The same inductive proof shows that if $A(x)$ and $W(x;t)$ have their
coefficients in a subfield $\cf$ of $\cs$, then the $P_{i,k}$ have their
coefficients in $\cf(\alpha_1, \ldots , \alpha _d)$.
\cqfd

\subsubsection{Proof of Theorem~\ref{General-endpoint}}
We are now ready to prove the main result of this section.
Let $m$ be the smallest integer (positive or not) such that $A(x):=x^m
A_1(x)$ is a polynomial in $x$. Then $A(0)\not = 0$. 
By Lemma~\ref{fixed-ordinate-general},
$$S_j(x;t)=\frac{1}{(2tA(x))^j}
\left( \frac{x^{mj}f_j^+(x;t)}{\sqrt{\Delta(x)}}
-
{{x^{mj}f^-_j(x;t)}{\sqrt{D\bar\Delta(\bar x)}}}\right).$$
 The positive part of the expression between brackets
is of the form
$$\frac{x^{mj}f^+_j(x;t)}{\sqrt{\Delta(x)}}-h_j(x;t),$$
where $h_j(x;t)$ is a Laurent polynomial in $x$ with coefficients in
the extension of $\qs(t)$ generated by $\sqrt D$ and the coefficients
of $\Delta(x)$ and $\bar \Delta(\bx)$. We then apply
Proposition~\ref{extract-positive} to obtain the expression of
$S_j^+(x;t)$. The value of $S_j^-(x;t)$ follows from the fact that
$S_j=S_j^-+S_j^+$. 

Finally, let us prove that $S_{i,j}$ belongs to $\qs _A$. 
With the above definition of $m$, 
$$S_j^+(x;t)=\frac{1}{(2tA(x))^j}
\left( \frac{x^{mj}f_j^+(x;t)}{\sqrt{\Delta(x)}}
-x^{mj}g_j(x;t)\right),$$
where $A(x)$ is a polynomial in $x$ such that  $A(0)\not = 0$. 
The part of the series that lies between the brackets can be seen,
either as a series in $t$ with coefficients in $\qs[x]$, or as a
series in $x$ with coefficients in $\qs _A$. Differentiating this
expression of $S_j^+(x;t)$ with respect to $x$, and setting $x=0$,
proves that for $i \ge 0$,  $S_{i,j}$ belongs to $\qs _A$. A similar
argument, starting from  $S_j^-(x;t)$, gives the announced result for
$i<0$. 
\cqfd
The extraction procedure will be applied explicitly to $S_1$ and
$S_2$, for the diagonal square lattice, in Section~\ref{autravail}.
We have implemented it as a \texttt{Maple} program and tested it on 
different models.

\subsubsection{Proof of Theorem~\ref{dim0}} 
\label{carre-sij}
Let us apply Theorem~\ref{General-endpoint} to the ordinary square
lattice.
As already observed, $A_1(x)=1$, so that we are not bothered by roots
of $A_1(x)$. Moreover, 
the symmetry of the model with respect to the line $x=0$ implies that
$\Delta(x) =\bar\Delta(x)$. Therefore the field
$\qs _A$ is generated by  $\sqrt D$ and the coefficients
of $\Delta(x)$. Eqs.~\Ref{Delta-def}, \Ref{delta-Delta} and
\Ref{Catalan-u} give
$$  \Delta(x)=
\left( 1-xu\ \frac{1+u}{1-u}\right) 
\left( 1-xu\ \frac{1-u}{1+u}\right) 
\quad \hbox{ and } \quad 
  \sqrt D
=t/u, $$
so that $\qs _A=\qs(u)$. Theorem~\ref{General-endpoint} states
that $S_{i,j}(t)$ belongs to $\qs(u)$. 
  In particular, expanding  in $x$ 
(or $\bx$) the series $S_0(x;t)=1/\sqrt{\Delta(x)}$, $S_1^+(x;t)$ and
$S_1^-(x;t)$ (given by~\Ref{S1+carre} and~\Ref{S1-carre}),
provides the announced 
expressions of $S_{1,0}$, $S_{0,1}$,  $S_{1,1}$ and  $S_{-1,1}$.

As any element of $\qs(u)$, the series  $S_{i,j}$ is either rational,
or quadratic, or quartic. Let us rule out rationality. Take a walk $w$ going from the point $(1,0)$ to $(i,j)$, and avoiding the
forbidden half-line. Let $m$ be the length of $w$. By adding $w$ at
the end of any walk ending at $(1,0)$, we obtain that for all $n\ge
m$,
$$a_{1,0}(n-m)\le a_{i,j}(n).$$
Similarly, by reversing the direction of $w$,
$$ a_{i,j}(n)\le a_{1,0}(n+m).$$
But  $a_{1,0}(n)=C_n\sim c\, 4^n n^{-3/2}$ (for $n$ odd). Hence there
exist two positive constants $c_1$ and $c_2$ such that, assuming
$n=i+j$ mod $2$,
$$ c_1 4^n n^{-3/2}\le a_{i,j}(n)\le  c_2 4^n n^{-3/2}.$$
 This cannot be the asymptotic behaviour of the coefficients
of a rational series.
\cqfd

\noindent{\bf Remark.} When we take into account the number of
vertical steps in the enumeration,  like in Section~\ref{section-refined}, some
expressions become more intricate (for instance, the series
$S_{0,1}(t,v)$ is of degree 8 while it is quadratic when $v=1$),
but certain results remain simple.
In particular, one derives directly from the fact that $S_0(x;t,v)=
\Delta(x)^{-1/2}$, where $\Delta(x)$ is
given in~\Ref{Delta-v} that
\begin{eqnarray*}
S_{1,0}(t,v)&=&\frac12 \left( X_3^{-1}+ X_4^{-1}\right) 
=\frac12 \left( X_1+X_2\right) \\
&=&   \frac{2-\sqrt{(1-2t-2tv)(1+2t-2tv)}
-\sqrt{(1-2t+2tv)(1+2t+2tv)}}{4t}.
\end{eqnarray*}
The series $X_1$ satisfies $X_1=t(1+2vX_1+X_1^2)$. Hence its
coefficients can be easily computed by the Lagrange inversion
formula. Moreover, $X_2(t,v)=X_1(t,-v)$. Finally,   the number of walks of length $2n+1$
going from $(0,0)$ to $(1,0)$ and having $2k$ vertical steps is found
to be
$$
{2n\choose 2k}2^{2k}C_{n-k},$$ in
accordance with the identity
$$
\sum_{k=0}^n{2n\choose 2k}2^{2k}
C_{n-k}=C_{2n+1}=a_{1,0}(2n+1).$$ This 
result suggests the existence
of a bijection between our walks and {\em bicolored Motzkin walks\/}
\cite[Ex.~2.2]{pizaler},
that would decrease the length by one and take the number of {\em
vertical\/} steps of our walk to the number of {\em horizontal\/}
steps of the Motzkin walk. Such a bijection has recently been described
in~\cite{florence}.

\subsubsection{Proof of Theorem~\ref{dim0bis}} 
\label{autravail}
For the diagonal square lattice, the canonical factorization of
$\delta(x)$ is given by~\Ref{Delta-bis}
and~\Ref{delta-Delta-2}. Again, $\Delta(x)=\bar \Delta(x)$.  These
values imply that the extension of $\qs(t)$ generated by $\sqrt D$ and
the coefficients of $\Delta(x)$ is simply
$\qs(t,C(4t^2))=\qs(t,\sqrt{1-16 t^2})$, ans has degree $2$ over
$\qs(t)$. However, $A_1(x)=x+\bx$ has two roots $\pm i$, so that the
field $\qs_A$ that contains all the series $S_{i,j}$ also contains
$\sqrt{\Delta(\pm i)}$ (and $i$). Using~\Ref{Catalan-u}, the series
$\Delta(x)$ and $D$ can be written in terms of $u$: 
$$
\Delta(x)=\bar\Delta(x)
=1-\frac{4u^2x^2}{(1-u^2)^2}
\ \ \ \hbox{ and }\ \ \ \sqrt D 
=\frac{(1-u^2)^2}{(1+u^2)^2}.
$$
In particular,  $\sqrt{\Delta(\pm i)}=(1+u^2)/(1-u^2)$. This series is of
degree $4$,  so that finally
$\qs_A=\qs(i,u)$. Theorem~\ref{General-endpoint} implies that
$S_{i,j}$ belongs to 
$\qs(i,u)$. But a series with {\em real\/} coefficients
belonging to $\qs(i,u)$    also belongs to $\qs(u)$. This proves
the first part of Theorem~\ref{dim0bis}. Rationality is ruled out as
for the ordinary square lattice. 

Let us now derive the announced expressions of the series $\smash{S_{i,j}}$.
The value of $S_{2i,0}$ follows from the fact that
$S_0(x;t)=1/\sqrt{\Delta(x)}$. The other three series require to apply
explicitly the extraction procedure described above, and actually
provide a good illustration of it.
 The series $S_1(x;t)$ is
given by~\Ref{S1-diagonal}. Observe that, for this model, $S_1(x;t)$
has only odd powers of $x$, and in particular, no constant term. Let
us divide 
it by $x$, in order to have even series,
$$
\frac{S_1^+(x;t)}{x} = \left( \frac{S_1(x;t)}{x}\right) ^+ 
=  \left(\frac{U(x;t)}{2t(1+x^2)}\right) ^+ 
$$
where
$$U(x;t)=\frac{1}{\sqrt{\Delta(x)}} -\sqrt{D\Delta(\bx)}.$$ This
remark allows us to work with the variable $X=x^2$, and to apply 
only once
 Corollary~\ref{diviseur-simple}.  The positive part of $U(x;t)$ is
$U^+(x;t)=\Delta(x)^{-1/2}-D^{1/2}$; with the notations of
Corollary~\ref{diviseur-simple}, $\alpha = -1$, $m=1$ and we obtain
$$
\frac{S_1^+(x;t)}{x} 
= \frac{1}{2t(1+x^2)}\left(\frac{1}{\sqrt{\Delta(x)}}
-\frac{1}{\sqrt{\Delta(i)}}\right),
$$
that is
$$S_1^+(x;t)= \frac{x}{2t(1+x^2)}\left(
\frac{1}{\sqrt{\Delta(x)}}-\frac{1-u^2}{1+u^2}\right).$$
By difference, this gives:
$$S_1^-(x;t)= \frac{\bx}{2t(1+\bx^2)} \left(\frac{1-u^2}{1+u^2}
-\sqrt{D \Delta(\bx)}\right).$$
The values of $S_{1,1}(t)$ and  $S_{-1,1}(t)$ follow, by extracting
the coefficient of $x$ (or $\bx$) in these series.

Finally, we need to evaluate the positive part of $S_2(x;t)$ in order
to compute $S_{0,2}(t)$. By Lemma~\ref{fixed-ordinate-general},
$$S_2(x;t)= \frac{U(x;t)}{2t^2(1+x^2) ^2},$$
where
$$U(x;t)= \frac{x^2-2t^2(1+x^2)^2}{\sqrt{\Delta(x)}}-x^2\sqrt{D \Delta(\bx)}.$$
Our first step is the extraction of the positive part of $U$:
$$U^+(x;t)= \frac{x^2-2t^2(1+x^2)^2}{\sqrt{\Delta(x)}}-x^2\sqrt{D}
+\frac{2u^2}{(1-u^2)^2}.$$
We now apply Corollary~\ref{diviseur-simple} (again, for the sake of efficiency, to
the variable $X=x^2$) with $\alpha = -1$ and  $m=2$. We need:
$$U^+(i;t)=\frac{-1}{\sqrt{\Delta(i)}}+\sqrt{D}
+\frac{2u^2}{(1-u^2)^2}$$
and
$$\frac{\partial U^+}{\partial
X}(i;t)=\frac{1}{\sqrt{\Delta(i)}}-\frac{2u^2}{(1-u^2)^2\Delta(i)^{3/2}}-\sqrt{D}.$$
Corollary~\ref{diviseur-simple} gives:
\begin{eqnarray*}
S_2^+(x;t)&=& \frac{1}{2t^2(1+x^2) ^2}\left(U^+(x;t)-U^+(i;t)-(x^2+1)
\frac{\partial U^+}{\partial X}(i;t)\right)\\
&=&  \frac{1}{2t^2(1+x^2) ^2}\left(\frac{x^2-2t^2(1+x^2)^2}{\sqrt{\Delta(x)}}
+\frac{2u^2(1-u^2)}{(1+u^2)^3} -x^2\frac{(1-u^2)(1+u^4)}{(1+u^2)^3}
\right).
\end{eqnarray*}
Taking the constant term in $x$ gives the announced expression of $S_{0,2}$.

\section{Other starting points}
\label{section-starting}
 From the results of Sections~\ref{complete} and \ref{prescribed},
some elementary operations 
on  walks allow us to solve a number of related enumerative
questions. We focus in this section on walks starting on the line
$y=0$ (but not necessarily at the origin), and avoiding the half-line
$\cal H$. The results we obtain will be useful in
Sections~\ref{kenyon-like} and 
\ref{asymptotics}, where probabilistic results are derived. Even though the
techniques we use here work for all models satisfying the conditions
of symmetry and small height variations, the results are easier to
state if we assume, in addition, that the set $A$ of steps is
symmetric with respect to reversal of the walks: $(i,j ) \in A
\Leftrightarrow (-i,-j) \in A$. This assumption is satisfied by the
ordinary and diagonal square lattices. Under this assumption, $A_0(x)
=A_0(\bx)$, and similarly $A_1(x)$ and $\delta(x)$ are symmetric in
$x$ and $\bx$. Consequently,
the canonical factorization of $\delta$ satisfies $\Delta(x) =
\bar\Delta(x)$.

We generalize the notations used so far by denoting $a_{i,j}^{[k]}(n)$
the number of $n$-step walks that start from $(k,0)$, end at $(i,j)$,
and have no contact with $\cal H$ except at their starting point (if
$k \le 0$). Similarly, we denote by $S^{[k]}(x,y;t)$ the corresponding
\gf :
$$S^{[k]}(x,y;t) = \sum_{n  \ge 0}\sum_{(i,j) \in \zs ^2} a_{i,j}^{[k]}(n)
x^i y^j t^n.$$ Section~\ref{negative} is devoted to the case $k\le 0$,
and Section~\ref{section-positive}  to the case $k > 0$.

\subsection{Starting on the ``forbidden'' half-line}
\label{negative}
\begin{propo}
\label{k<0}
With the notations of Section{\em ~\ref{complete}},
the \gf \ for walks that start on $\cal H$, but never return to this
half-line is:
$$S^-(z;x,y;t):= \sum _{k \ge 0} S^{[-k]}(x,y;t)z^k= 
\frac{\sqrt{D\Delta(\bx) \Delta(z)}}{(1-z\bx)K(x,y)}=
\frac{\sqrt{ \Delta(z)}}{1-z\bx}S(x,y;t).$$
\end{propo}
{\bf Proof.}
First, observe that $S^-(0;x,y;t)$ counts  walks starting at the
origin, and hence, coincides with the complete \gf \ $S(x,y;t)$ for
walks of the slit plane.
Now take a walk starting from
$(-k,0)$, with $k >0$, and move it one step to the right: we obtain a
walk of the same type starting from $(-k+1,0)$. Conversely, a walk
starting from  $(-k+1,0)$, translated one step to the left, provides
either a walk of the right type starting from $(-k,0)$, or the
concatenation $uvw$ of a (reversed) bridge 
 $u$  going from  $(0,0)$ to $(-k,0)$,
a (possibly empty) sequence $v$ of bridges ending at $(0,0)$, and a
usual walk $w$ on the slit plane (recall the definition of bridges, at
the beginning of Section~\ref{section-complete-general}). In terms of \gfs :
$$z\bx S^-(z;x,y;t)= S^-(z;x,y;t)-S^-(0;x,y;t)+
\frac{B(z)-B(0)}{1-B(0)}S(x,y;t).$$
We conclude thanks to Theorem~\ref{dim2-general} and
Proposition~\ref{bridges0}.
\cqfd
The following corollary will be useful in Section~\ref{kenyon-like}.
\begin{coro}\label{autrespoints}
 For  $(i,j) \in \zs ^2\setminus {\cal H}$, the \gf \  
 for walks that start at $(i,j)$,
end on $\cal H$, but otherwise avoid $\cal H$, is
$$
 D_{i,j}(z;t):= \sum_{n  \ge 0}\sum_{k \ge 0} a_{i,j}^{[-k]}(n)
z^k t^n =\sqrt{\Delta(z)} S^+_{i,j}(z;t), 
$$
where $S^+_{ i,j}(x;t)$ is the following  section of the complete \gf
\ $S(x,y;t)$:
$$S^+_{ i,j}(x;t)=\sum_{k\ge 0}x^{k}S_{i+k,j}(t).$$ 
\end{coro}
{\bf Proof.} By the symmetry assumption, the series $D_{i,j}(z;t)$ is
the coefficient of $x^i 
y^j$ in $S^-(z;x,y;t)$. The result
follows from the last expression of  $S^-(z;x,y;t)$ given in Proposition~\ref{k<0}.
\cqfd

\subsection{Starting at $(k,0)$, with $k>0$}
\label{section-positive}
\begin{propo}
\label{k>0}
With the notations of Section{\em ~\ref{complete}},
the \gf \ for walks that start on the half-line $\{(k,0), k>0\}$ and
avoid $\cal H$ is:
$$S^+(z;x,y;t):= \sum _{k > 0} S^{[k]}(x,y;t)z^k= 
\frac{zx}{(1-zx)K(x,y)}\sqrt{\frac{\Delta(\bx)}{\Delta(z)}}=
\frac{zx S_0(z;t)S(x,y;t)}{(1-zx)\sqrt{D(t)}}.$$
\end{propo}
{\bf Proof.} We essentially copy the argument we used in the proof of
Proposition~\ref{k<0}. By moving one step to the left a walk counted by
$S^+(z;x,y;t)$, we obtain either a walk of the same type, or the
concatenation $uvw$ of a (reversed) walk on the slit plane ending at
$(k,0)$, a possibly empty sequence $v$ of bridges ending at
$(0,0)$, and a walk $w$ on the slit plane. In terms of \gfs ,
$$\bx \bz S^+(z;x,y;t)= S^+(z;x,y;t)+
\frac{S_0(z;t)S(x,y;t)}{1-B(0)}.$$ 
The result follows, using  Theorem~\ref{dim2-general} and
Proposition~\ref{bridges0}.
\cqfd

Let us call {\em loop\/} a walk that starts and ends at 
the same point of the positive $x$-axis and avoids $\cal H$. Let $L_k(t)$ be
the length \gf \ for loops starting and ending at $(k,0)$.
\begin{coro} \label{loops} The \gf \ for loops is
$$L(z;t):=\sum_{k >0} L_k(t) z^k= \frac{z}{(1-z)\sqrt {D(t)}}\sum_{k \ge 0}
S_{k,0}(t)^2 z^k.$$
\end{coro}
{\bf Proof.} The series $L_k(t)$ is the
coefficient of $z^k x^k y^0$ in the series $S^+(z;x,y;t)$. Let us
consider the last expression of this series given in
Proposition~\ref{k>0}: the coefficient of $y^0$
in $S(x,y;t)$ being $S_0(x;t)$, the result easily follows.
\cqfd
Our last proposition will be used in Section~\ref{asymptotics}
 to obtain the average number
of visits of a (long) walk on the slit plane to the point $(k,0)$.
\begin{propo}
\label{visites}
For $k>0$, the \gf \ for walks on the slit plane (starting from
$(0,0)$) that visit the 
point $(k,0)$  is
$$ V_k(x,y;t,v)=\frac{v S_{k,0}(t)
S^{[k]}(x,y;t)}{L_k(t)^2\left(1-v(1-1/L_k(t))\right)}.$$ 
This series counts walks by their length (variable $t$), number of visits to
$(k,0)$  (variable $v$), and position of their endpoint (variables $x,y$).
\end{propo}
Recall that $S_{k,0}$ counts walks on the slit plane ending at $(k,0)$,
$S^{[k]}(x,y;t)$ counts walks starting from  $(k,0)$, and $L_k(t)$
counts loops  starting from  $(k,0)$. These series can be computed
respectively 
from the expansions of $S_0(x;t)$ (Theorem{~\ref{dim2-general}}),
of $S^+(z;x,y;t)$ (Proposition{~\ref{k>0}}),  and of $L(z;t)$
(Corollary{~\ref{loops}}). 

\medskip
\noindent {\bf Proof.}
Let us say that a non-empty loop is {\em primitive\/} if it visits
exactly twice its starting (and ending) point. Let $P_k(t)$ be the
length \gf \ for primitive loops starting at $(k,0)$. Clearly,
$L_k=(1-P_k)^{-1}$. Equivalently, $P_k =1-1/L_k$.

Similarly,  a walk on the slit plane ending at $(k,0)$
is {\em primitive\/} if it visits only once its endpoint. The \gf \ for
primitive walks ending at  $(k,0)$ is 
 $S_{k,0}/L_k$. Finally,
$S^{[k]}/L_k$ counts walks starting from $(k,0)$ that never return to
their starting point.

 A walk on the slitplane visiting $(k,0)$ can be seen in a unique
way as the concatenation $up_1 p_2 \cdots p_\ell w$ of a primitive walk
$u$ going from the origin to $(k,0)$, a sequence of $\ell$ primitive
loops starting at $(k,0)$, and a primitive walk $w$ starting from $(k,0)$.
Such a walk visits exactly $\ell+1$ times the point $(k,0)$. The
result follows.
\cqfd

\section{The hitting distribution of a half-line} \label{kenyon-like}
We focus in this section on the ordinary square lattice.
The  results obtained in Sections~\ref{prescribed}  and \ref{section-starting}
allow us to solve a
number of probabilistic 
questions ``\`a la Kenyon''.  Let $(i,j)$ be a point of $\zs^2$.
A random walk starting from
$(i,j)$ hits the half-line ${\cal H}=\{(k,0), k\le 0\}$ with probability
$1$. The probability that 
the {\em first\/} hitting point is $(0,0)$ is
\beq p_{i,j} =\sum_{n\ge 0} \frac{a_{i,j}(n) }{4^n} =S_{i,j}(1/4).
\label{pij-kenyon}\eeq
More precisely, $a_{i,j}(n)/4^n$ is the probability that this event occurs
after $n$ steps. Theorem~\ref{dim0} states that  $S_{i,j}(t)$ belongs to
$\qs(u)$, where the series $u$ is given by~\Ref{u-def}. As $u=\sqrt 2
-1$ when $t=1/4$, this theorem implies 
that $p_{i,j} \in \qs[\sqrt{2}]$ and
in particular, that $p_{0,1}=1/2$ and $p_{1,0}=2-\sqrt{2}$, as
stated in R. Kenyon's e-mail. 
We have written a program to compute $p_{i,j}$:
for $|i|+|j| \le 10$,  the probability   $p_{i,j}$ is irrational,
unless $(i,j)=(0,1)$ or $(0,-1)$.

Lemma~\ref{fixed-ordinate-general} and Theorem~\ref{dim1} tell us
how to compute the series
$$S_j^+(x;t)= \sum_{i \ge 0} S_{i,j}(t) x^i \quad \hbox{ and } \quad
S_j^-(x;t)= \sum_{i < 0} S_{i,j}(t) x^i.$$
Setting $t=1/4$ in these expressions  
provides  explicit values  for
\gfs \ of the form
$$\sum_{i \ge 0}p_{i,j}x^i \ \ \ \hbox{ and }\ \ \ \sum_{i <
0}p_{i,j}x^i .$$
For instance,  $S_0(x)=
{\Delta(x)}^{-1/2}$, where $\Delta(x)$ is given by~\Ref{Delta-def}, and consequently,
$$\sum_{i \ge 0}p_{i,0}x^i= \frac{1}{\sqrt{(1-x)\left(1-x(\sqrt
2-1)^2\right)}}.$$ 

More generally, given $k\ge 0$,
one can  ask about the probability $p_{i,j}^{[k]}$
that the first hitting point of a random walk starting from $(i,j)$ is
$(-k,0)$. For $k$ fixed, the function 
$f$ defined by $f(i,j)=p_{i,j}^{[k]}$ is the unique bounded function
on $\zs ^2$, harmonic  on $\zs ^2\setminus {\cal H}$, such that
$f(-\ell, 0)=\delta_{k,\ell}$  for $\ell \ge 0$
(see~\cite[Theorem~1.4.8]{lawler}).
 The {\em hitting
distribution\/} of the half-line $\cal H$, starting from $(i,j)$, is
condensed in the following series,
$$ \sum_{k\ge 0}p_{i,j}^{[k]}z^k=D_{i,j}(z;1/4),$$
where we have used the notations of Corollary \ref{autrespoints}. This
corollary, combined with Theorem~\ref{dim1}, allows us to compute
this series explicitly, for a fixed value of $(i,j)$.
\begin{Theorem}
\label{hitting-theorem}
Let $(i,j)\in \zs^2\setminus {\cal H}$. The hitting
distribution  of the half-line $ {\cal H}$, starting from
$(i,j)$, is of the following form:
$$ \sum_{k\ge 0} p_{i,j}^{[k]}z^k 
=f(z)-g(z) \sqrt{(1-z) \left ( 1-z( \sqrt{2}-1) ^2\right) }
$$
where $f(z)$ and $g(z)$ are Laurent polynomials in $z$, with
coefficients in $\qs$ and $\qs[\sqrt 2]$ 
respectively, that satisfy $f(1)=1$, $g(1)
\not = 0$.  Consequently, as $k\rightarrow \infty$, the probability
that the first hitting point is $(-k,0)$ is
$$p_{i,j}^{[k]}\sim  g(1) \sqrt{\frac{\sqrt 2-1}{2 \pi}}\ k^{-3/2}$$
and
 the probability that the
first hitting point is to the left of $(-k,0)$ is
$$\sum_{\ell \ge k} p_{i,j}^{[\ell]}\sim  g(1) \sqrt{\frac{2(\sqrt
2-1)}{\pi}}\ k^{-1/2}.$$ 
The Laurent polynomials $f(z)$ and $g(z)$  can be computed
explicitly. For instance, $f(z)=g(z)=\bz$ if 
$(i,j)=(1,0)$, so that in this case $g(1)=1$. 
\end{Theorem}
{\bf Proof.}
Corollary~\ref{autrespoints} expresses $D_{i,j}(z;t)$ in terms of the section
$S_{i,j}^+(z;t)$. The latter series only differs from the section
$S_j^+(z;t)/z^i$ (the generic form of which is given by
Theorem~\ref{dim1}) by a finite number of series $S_{k,j}$; in other
words, according to Theorem~\ref{dim0}, by a Laurent polynomial
$h(z;u)$ in $z$ with 
coefficients in $\qs(u)$. Hence the combination of Corollary~\ref{autrespoints}
and Theorem~\ref{dim1} imply that  $D_{i,j}(z;t)$ is of the following
form:
\beq D_{i,j}(z;t)=f(z;t) - g(z;u) \sqrt{\Delta(z)},
\label{dij-generic}\eeq
where $f$ is a  Laurent polynomial in $z$ and $t$ with rational
coefficients, and $g$ is a  Laurent polynomial in $z$ with
coefficients in $\qs(u)$. In particular, for $(i,j)=(1,0)$, we have
$S^+_{1,0}(z;t)=(S_0(z;t)-1)/z$ and 
$$  D_{1,0}(z;t) =
 \frac{1-\sqrt{\Delta(z)}}{z}. $$
Setting $t=1/4$ (and $u=\sqrt 2 -1$) in~\Ref{dij-generic} provides the
announced form for $D_{i,j}(z;1/4)= \sum_k p_{i,j}^{[k]}z^k$.

The recurrence of random walks on $\zs ^2$ implies that the half-line
is visited almost surely, so that
$$ \sum_{k\ge 0} p_{i,j}^{[k]}=f(1)=1.$$ 

 The asymptotic behaviour of $ p_{i,j}^{[k]}$ follows from the nature
of the singularities of the series  $ D_{i,j}(z;1/4)$. The dominant
singularity is at $z=1$, and is a square root singularity, provided
that $g(1)\not = 0$.  This is, at least, the case when $(i,j)=(1,0)$.

In general, if $ z=1$ was a root of $g(z)$, say, of multiplicity $m$, then
$p_{i,j}^{[k]}$ would decay like $k^{-m-3/2}$, that is, much faster than
$p_{1,0}^{[k]}$.  As in Section~\ref{carre-sij}, considering a (fixed) walk
going from $(i,j)$ to $(1,0)$ on the slit plane proves that this is
impossible. 

 Similarly, the study of the behaviour of  $(1-D_{i,j}(z;1/4))/(1-z)$
around its dominant
singularity provides the asymptotic
behaviour of the probability
that the hitting abscissa is smaller than $-k$.
\cqfd

\bigskip

\noindent {\bf Remarks}  

\noindent 1. This result has to be compared with Lemma~$6$
in \cite{kesten}, where it is proved (with our notations), that for $0
\le i \le k/2$, 
\beq \sum_{\ell \ge k} p_{-i,1}^{[\ell]}\le c\ ((i+1)k)^{-1/2},
\label{bounds} \eeq
for a constant $c$ independent 
of $i$ and $k$. This statement
contains some uniformity in $i$ which is absent from our result. It also
implies that, for all $i$ and $j$, there exists a constant $c_{i,j}$
such that for all $k \ge 0$,
$$\sum_{\ell \ge k} p_{i,j}^{[\ell]}\le c_{i,j}\ k^{-1/2},$$
but this, in turn, is weaker than Theorem~\ref{hitting-theorem}.
Related results are described in \cite[Section~2.4]{lawler}. For
instance, the bound~\Ref{bounds} follows directly from Eq.~(2.40) in
\cite{lawler}.
 Analogously, the probability that a planar brownian
motion starting from $(1,0)$  hits the half-line for the first time at
abscissa smaller than $-k$ also decays like $k^{-1/2}$
(see e.g. \cite[with $\alpha=2$]{ray}).

\noindent 2. The probability $p_{i,j}^{[k]}$ admits another probabilistic
interpretation. Let $(W_n)_{n \ge 0}$, be a random walk on the square
lattice starting from $W_0=(-k,0)$. Let $T=\min\{n >0 : W_n\in{\cal
H}\}$, and let $V$ denote the number of visits of the walk to the
point $(i,j)$ before $T$. Then the expectation of $V$ is
\begin{eqnarray*}
\GE(V)&=& \sum_{n >0} \Pr(n< T, W_n=(i,j))\\
&=&  \sum_{n >0}\frac{1}{4^n}\,a_{i,j}^{[-k]}(n) \\
&=&p_{i,j}^{[k]}.
\end{eqnarray*}
As $p_{i,j}^{[k]}$ has been seen to be a probability,  this average
number of visits is always less than $1$.

\section{Properties of long walks on the slit plane}
\label{asymptotics}
Again, we focus on the ordinary square lattice.
\subsection{Transience and Green function}\label{section-transient} 
It is well-known that  random walks on the square lattice are
{\em recurrent\/}: any given point $(i,j)$ of the lattice is visited
 with probability $1$, and is actually visited infinitely many
times. In more enumerative terms, the 
proportion of walks of length $n$ visiting $(i,j)$ tends to $1$ as $n$
goes to infinity, and the average number of visits of $n$-step walks
to this point tends to infinity. 

This is no longer the case for walks on the slit plane. The forbidden
half-line creates a long-range repulsion of the walks, which become
transient. This result is not surprising, 
and can probably be proved by various methods; but the enumerative results
we have obtained, and more 
especially Proposition~\ref{visites}, can be used to obtain
{\em exact quantitative information\/}. For instance, we can compute,
for any point 
$(k,0)$ with $k>0$, the probability that it is visited by an $n$-step
walk. This probability is
\beq \frac{[t^n]V_k(1,1;t,1)}{a(n)}= \frac{1}{a(n)}[t^n]
\left( \frac{S_{k,0}(t)S^{[k]}(1,1;t)}{L_k(t)}\right).
\label{pk}\eeq
In this expression, $a(n)$ denotes the total number of $n$-step walks
on the slit plane, 
and the notation $[t^n]$ means ``the coefficient of $t^n$''. 
We shall prove
that this probability converges as $n$ goes to infinity to a limit
that is strictly less than $1$.
By differentiating $V_k$ with respect to $v$, we can also compute the
average number of visits to $(k,0)$ of $n$-step walks, which is
\beq \frac{[t^n]V_k'(1,1;t,1)}{a(n)}= \frac{1}{a(n)}[t^n]
\left( {S_{k,0}(t)S^{[k]}(1,1;t)}\right).\label{vk}\eeq
Again, this quantity 
will be shown to converge to a {\em finite\/}
limit as $n$ goes to infinity. By analogy with ordinary random walks,
we call it the value at $(k,0)$ of the {\em Green function\/} of our
model.\\
Recall that for any $i,j$, the series $S_{i,j}(1/4)$ is always finite
(and at most $1$, by~\Ref{pij-kenyon}). 
\begin{propo}
 As $n$ tends to infinity, the proportion of walks of length $n$
 visiting $(k,0)$ tends to
$$p_k= S_{k,0}(1/4) \ \frac{\sum_{\ell =0}^{k-1}S_{\ell,0}(1/4)}
{\sum_{\ell =0}^{k-1}S_{\ell,0}(1/4)^2} \ <1,$$
 while the average number of visits to this point goes to
$$v_k= 4(\sqrt2-1) \,S_{k,0}(1/4) \sum_{\ell =0}^{k-1}S_{\ell,0}(1/4).$$
These expressions, together with 
$$S_0(x;1/4)=\sum_{\ell \ge 0}S_{\ell,0}(1/4)x^\ell= \frac{1}{\sqrt{(1-x)\left(1-x(\sqrt 2-1)^2\right)}},$$
allow us to compute $p_k$ and $v_k$ for any $k>0$. For instance, walks
are more likely to visit $(1,0)$ than $(2,0)$:
$$p_1=2-\sqrt2\approx0.586 , \ \ \
p_2=\frac{5}{34}(19-11\sqrt2)\approx.506$$
but spend more time at  $(2,0)$ than  $(1,0)$:
$$v_1=4(3\sqrt2-4)\approx0.97, \ \ \ 
v_2=10(22\sqrt2-31)\approx1.13.$$
\end{propo}
{\bf Proof.}
 From the last expression of $S^+(z;x,y;t)$ given in
Proposition~\ref{k>0}, we obtain:
$$S^{[k]}(1,1;t)=\frac{S(1,1;t)}{\sqrt {D(t)}} \sum_{\ell
=0}^{k-1}S_{\ell,0}(t).$$ 
Similarly, Corollary~\ref{loops} gives
$$L_k(t)=\frac{1}{\sqrt {D(t)}} \sum_{\ell
=0}^{k-1}S_{\ell,0}(t)^2.$$ 
The series $D(t)$ is given by~\Ref{delta-Delta}, which implies
$1/\sqrt{D(1/4)}= 4(\sqrt2-1)$.
The announced expressions of $p_k$ and $v_k$  then
follow from~\Ref{pk} and~\Ref{vk} by analysis of the singularities of
the series into consideration.

Finally, one derives from the expression of $S_0(x;1/4)$ that
$S_{k,0}(1/4)$ is a strictly decreasing function of $k$. Consequently, 
$p_k<1$.
\cqfd

\subsection{Limit law for the  coordinates of the endpoint}
\label{limit-laws}
\begin{figure}
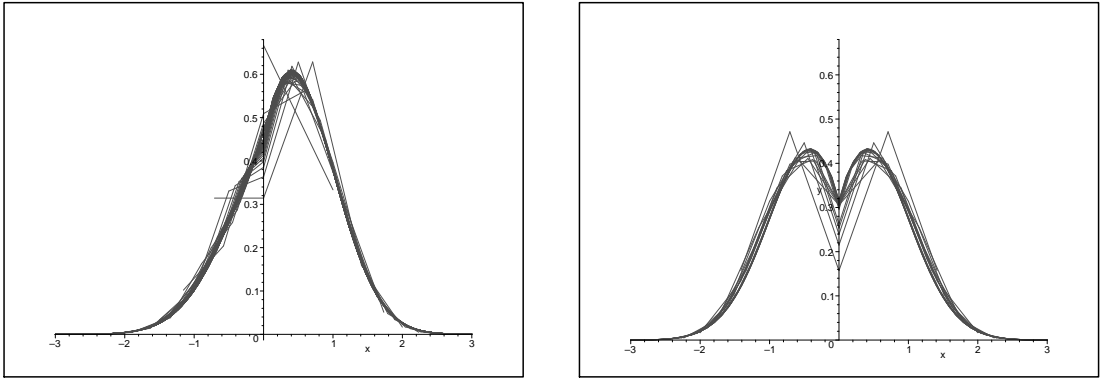

\begin{center}
\begin{turn}{-90}
\epsfxsize=5cm
\epsf{limiteX.eps}
\end{turn}
\hspace{5mm}
\begin{turn}{-90}
\epsfxsize=5cm
\epsf{limiteY.eps}
\end{turn}
\end{center}
\caption{The convergence of $\sqrt{n}\Pr(X_n=i)$, plotted against
$i/\sqrt n$ (left), and $\sqrt{n}\Pr(Y_n=j)$ against $j/\sqrt{n}$
(right), for $n=2,3, \ldots , 10, 20, 30, \ldots , 100$.
}
\label{limiteY}
\end{figure}

When all walks of length $n$ on the slit plane are
taken equally likely, the coordinates of their endpoints
become random variables $X_n$ and $Y_n$.
It is well-known that for an ordinary random walk,
these coordinates, normalized by $\sqrt n$, 
converge to a two-dimensional  centered normal law.
 This suggests to try the same 
 normalization for walks on the slit plane.
By expanding in $t$ the complete \gf \ of
Theorem~\ref{dim2}, the probabilities 
\[
\Pr((X_n,Y_n)=(i,j))\;=\;\frac{a_{i,j}(n)}{a(n)}
\]
 can be explicitly computed
for small values of $n$ and any $(i,j)$. 
The plots of the marginals,
$\sqrt n\Pr(X_n=i)$ against $i/\sqrt n$, and $\sqrt{n} \Pr(Y_n=j) $
against $j/\sqrt{n}$,  shown on Fig.~\ref{limiteY}, suggest that the
normalized random variables $X_n/\sqrt n$ and $Y_n/\sqrt{n}$ also
converge in distribution.
These plots actually even suggest the
existence of a {\em local\/} limit law. We have proved this for the ordinate
$Y_n/\sqrt n$, but not for the abscissa, and we shall simply prove
here the convergence in distribution of the normalized endpoint.
 Unsurprisingly, the limit law we obtain corresponds to a
two-dimensional Brownian conditioned (with care) not to hit a
half-line~\cite{legall}, and is related to the solution of the
associated diffusion equation (see~\cite[Eq.~(29)]{considine}). 

\begin{figure}
\begin{center}
\epsfxsize=5cm
\epsfysize=4cm
\epsf{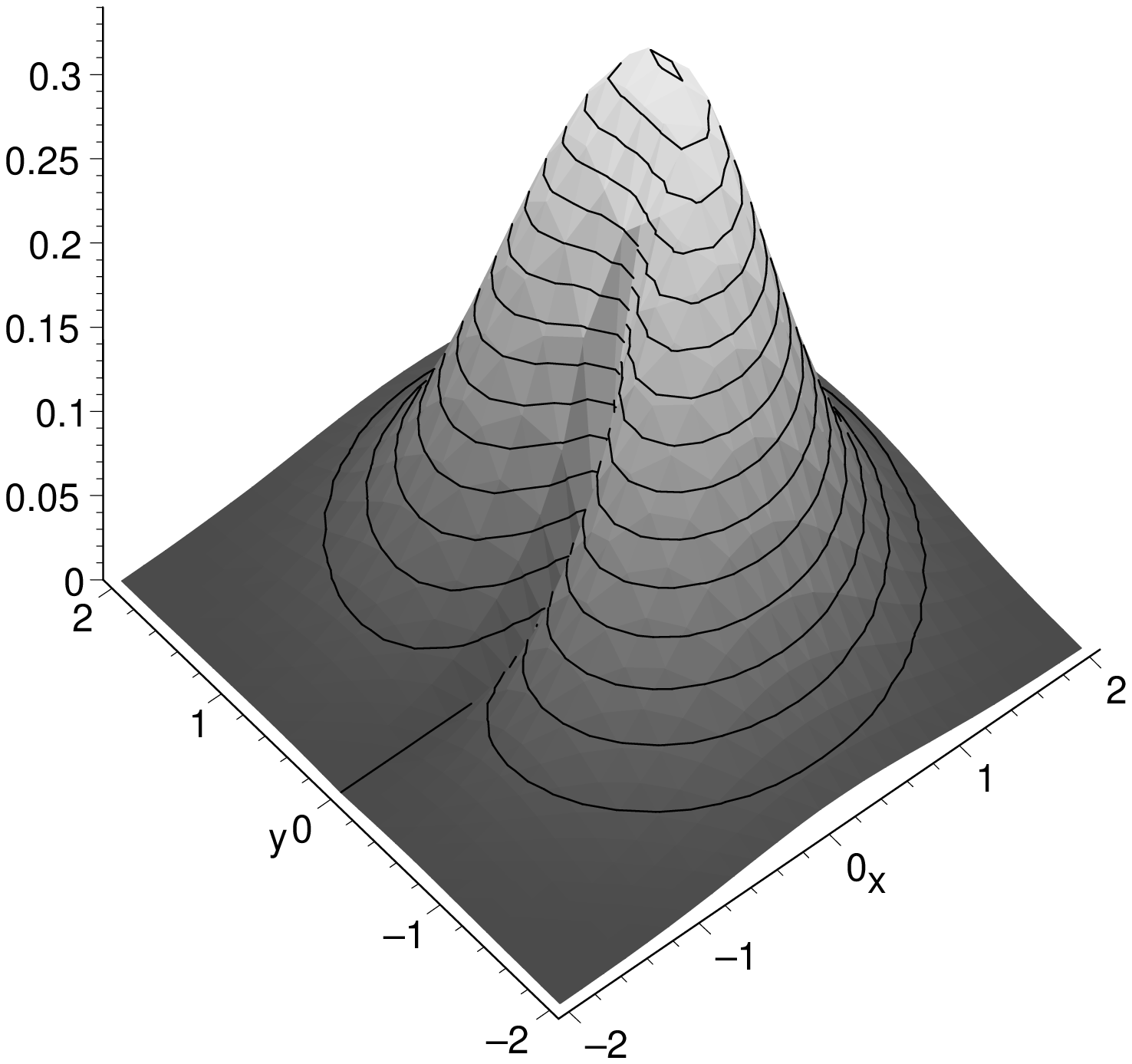}
\hspace{15mm}
\epsfxsize=5cm
\epsfysize=5cm
\epsf{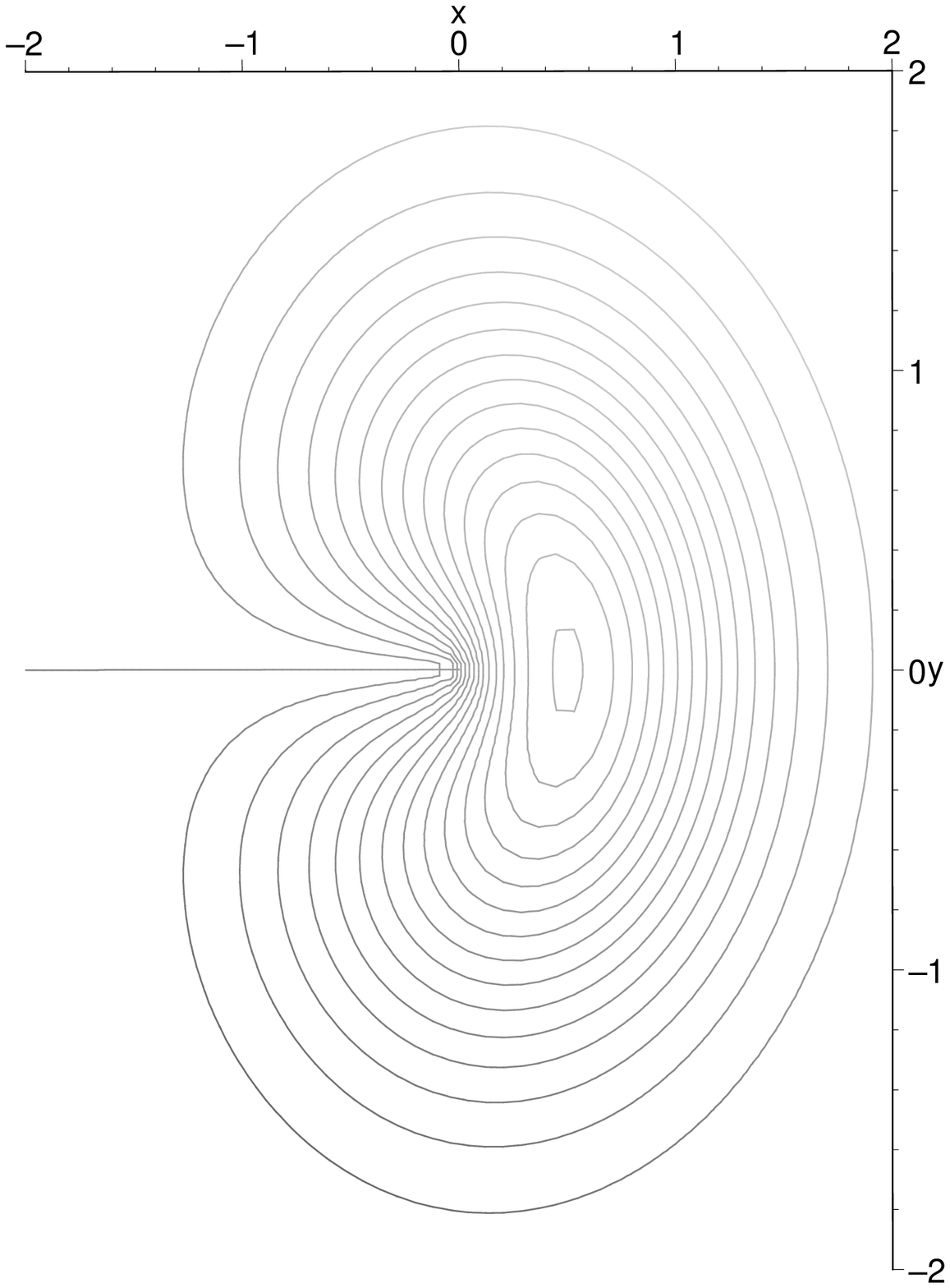}
\end{center}
\caption{The density $f(x,y)$ of $(X,Y)$.
\label{bivariate}}
\end{figure}

\begin{Theorem}
\label{distribXY}
 The sequence of joint random variables $(X_n/\sqrt n,Y_n/\sqrt
n)$ converges in distribution towards a pair $(X,Y)$ of density
$$
f(x,y)=\frac{\sqrt2}{\Gamma(1/4)}\;e^{-(x^2+y^2)}\sqrt{x+\sqrt{x^2+y^2}}
$$
with respect to the Lebesgue measure on $\rs ^2$. This density is
shown in Fig.{\em ~\ref{bivariate}}.  
\end{Theorem}
The pair $(X,Y)$
actually admits a simpler description in polar coordinates.
Let  $R=|X+iY|$ and $\Theta=\arg(X+iY) \in [-\pi,\pi [$.
Then $(R, \Theta)$ has  density
\[
g(\rho,\theta)=\frac{2}{\Gamma(1/4)}\;\rho^{3/2}e^{-\rho^2}\cos(\theta/2)
\]
with respect to the Lebesgue measure on $\mathbb{R}^+\times[-\pi,\pi]$.
Using this expression of the density, we can easily compute the moments of
$(X,Y)$. Thanks to Theorem~\ref{distribXY}, we obtain in particular the
following asymptotic results.
%
\begin{coro}
As $n$ goes to infinity, the endpoint $(X_n, Y_n)$  of a random
$n$-step walk on the slit plane satisfies:
$$
\GE( X_n)\sim\frac{\Gamma(3/4)}{\Gamma(1/4)}\sqrt n, 
\qquad
\GE(Y_n)=0, 
\qquad
\GE \left(\sqrt{X_n^2+Y_n^2}\right) \sim 
3\frac{\Gamma(3/4)}{\Gamma(1/4)}\sqrt
n,$$
$$\GE \left(X_n^2\right)\sim \frac7{12} n 
\qquad \hbox{and }\qquad
\GE \left(Y_n^2\right)\sim \frac23 n.$$
\end{coro}

\noindent {\bf Remark.} When the square lattice is replaced by the
diagonal lattice, elementary steps undergo a dilatation of a factor
$\sqrt2$. One can actually prove that Theorem~\ref{distribXY} holds for
the diagonal case, upon normalizing the variables $X_n$ and $Y_n$ by
$\sqrt{2n}$ instead of $\sqrt{n}$.

\begin{figure}[tb]
\begin{center}
\begin{picture}(0,0)%
\includegraphics{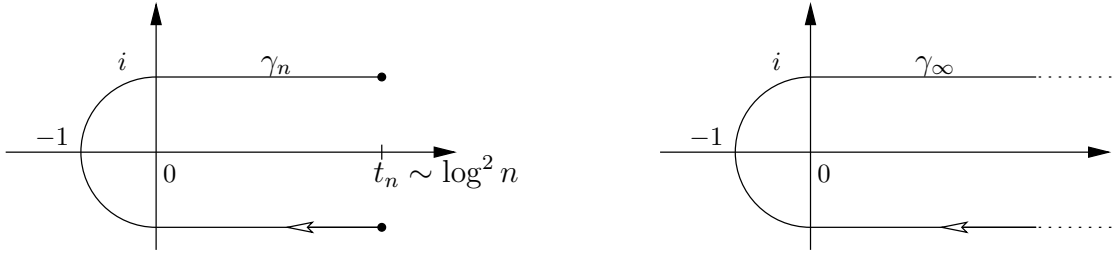}%
\end{picture}%
\setlength{\unitlength}{4144sp}%
\begingroup\makeatletter\ifx\SetFigFont\undefined%
\gdef\SetFigFont#1#2#3#4#5{%
  \reset@font\fontsize{#1}{#2pt}%
  \fontfamily{#3}\fontseries{#4}\fontshape{#5}%
  \selectfont}%
\fi\endgroup%
\begin{picture}(6639,1509)(799,-1558)
\put(5671,-1141){\makebox(0,0)[lb]{\smash{\SetFigFont{10}{12.0}{\familydefault}{\mddefault}{\updefault}
$0$
}}}
\put(5401,-466){\makebox(0,0)[lb]{\smash{\SetFigFont{10}{12.0}{\familydefault}{\mddefault}{\updefault}
$i$
}}}
\put(4906,-916){\makebox(0,0)[lb]{\smash{\SetFigFont{10}{12.0}{\familydefault}{\mddefault}{\updefault}
$-1$
}}}
\put(6256,-466){\makebox(0,0)[lb]{\smash{\SetFigFont{12}{14.4}{\familydefault}{\mddefault}{\updefault}
$\gamma_\infty$
}}}
\put(3016,-1141){\makebox(0,0)[lb]{\smash{\SetFigFont{12}{14.4}{\familydefault}{\mddefault}{\updefault}
$t_n\sim \log ^2 n$
}}}
\put(1756,-1141){\makebox(0,0)[lb]{\smash{\SetFigFont{10}{12.0}{\familydefault}{\mddefault}{\updefault}
$0$
}}}
\put(1486,-466){\makebox(0,0)[lb]{\smash{\SetFigFont{10}{12.0}{\familydefault}{\mddefault}{\updefault}
$i$
}}}
\put(991,-916){\makebox(0,0)[lb]{\smash{\SetFigFont{10}{12.0}{\familydefault}{\mddefault}{\updefault}
$-1$
}}}
\put(2341,-466){\makebox(0,0)[lb]{\smash{\SetFigFont{12}{14.4}{\familydefault}{\mddefault}{\updefault}
$\gamma_n$
}}}
\end{picture}
\end{center}
\caption{The contours $\gamma_n$ and $\gamma_\infty$.}
\label{gamma-t-infty}
\end{figure}

\medskip
Our proof of Theorem~\ref{distribXY} is (as the rest of this paper)
based on the exact enumerative results of Section~\ref{complete} and
provides an alternative 
to more classical probabilistic proofs. Let  $\Phi_n$ be the
characteristic function of $(X_n/\sqrt
n,Y_n/\sqrt n)$.  We shall see that the {\em analysis of
singularities\/} developed by Flajolet and
Odlyzko~\cite{flajolet-odlyzko} implies automatically  the
pointwise convergence of  $\Phi_n$ and gives an expression for its
limit.
 This limit will then be  identified as the characteristic
function of $(X,Y)$. 
\begin{propo}[Convergence of the characteristic functions]\label{charact}
Let $\sigma$ and $\tau$ be real numbers. Then 
\[
\lim_{n\rightarrow+\infty}\Phi_n(\sigma,\tau)=\Phi(\sigma,\tau),
\]
where 
$$
\Phi(\sigma,\tau)=
\frac{2}{i\Gamma(1/4)}
\int_{\gamma_\infty}
\frac{(2\sqrt{-t}+i\sigma)^{1/2}e^{-t}}{\sigma^2+\tau^2-4t}
dt,
$$
with $\gamma_\infty$ the Hankel contour around $[0,+\infty[$ 
shown on the right side of Fig.{\em~\ref{gamma-t-infty}}.
\end{propo}

\noindent\textbf{Proof.}
The characteristic function $\Phi_n(\sigma,\tau)$  can be expressed in
terms of the complete \gf \ $S(x,y;t)$:
$$ 
\Phi_n(\sigma,\tau)=E\left(e^{i(\sigma X_n+\tau Y_n)/\sqrt n}\right)
=\frac{[t^n]S(e^{i\sigma/\sqrt n},e^{i\tau/\sqrt n};t)}
{a(n)}
$$
where $a(n)$ is the total number of $n$-step walks on the slit plane.
 Recall that,
from Theorem~\ref{dim2}, this number grows like $4^n n^{-1/4}$, up to
 an explicit multiplicative 
constant. As we want to prove that $\Phi_n(\sigma,\tau)$ converges,
this means that the coefficient of $t^n$ in $S(e^{i\sigma /\sqrt
n},e^{i\tau /\sqrt n};t)$ also has to grow like  $4^n n^{-1/4}$. 

We shall estimate this coefficient thanks to the 
``analysis of singularities''~\cite{flajolet-odlyzko}. For the sake of
completeness, we give all the details of the calculation, but we
insist on the fact that it is a rather direct application
of~\cite{flajolet-odlyzko}.  We begin with
Cauchy's formula, and force the factor $4^n$ to appear by setting $t=z/4$:
\beq 
{2i\pi}\frac{a(n)}{4^n}\Phi_n(\sigma,\tau)=\int_{\cal C}
S(e^{i\sigma/\sqrt n},e^{i\tau/\sqrt n};z/4)\frac{dz}{z^{n+1}},
\label{cauchy}\eeq
where $\cal C$ is any simple contour positively encircling the origin,
inside the domain of analycity of $S(e^{i\sigma/\sqrt
n},e^{i\tau/\sqrt n};z/4)$. We expect this integral to behave like $n^{-1/4}$.
 By Theorem \ref{dim2},  
\beq 
S(x,y;z/4)
   =\frac{\left(2-z(1+e^{-i\sigma/\sqrt n})+2\sqrt{1-z}\right)^{1/2}
          \left(2+z(1-e^{-i\sigma/\sqrt n})+2\sqrt{1+z}\right)^{1/2}}
             {4-2z\left(\cos\left(\sigma/\sqrt n\right)
                       +\cos\left(\tau/\sqrt n\right)\right)}
\label{explicite}\eeq
where we  denote $x=e^{i\sigma/\sqrt n}$ and
$y=e^{i\tau/\sqrt n}$. 
We choose the principal determination of the square root on $\cs
\setminus \rs ^-$, given by
\[\sqrt{\rho
e^{i\theta}}=\sqrt{\rho}e^{i\theta/2}
\quad\textrm{ for }\;
\rho>0\;\textrm{ and }\; \theta\in \ ]-\pi,\pi[.
\]
The singularities of the series 
 $S(x,y;z/4)$ may, at first sight,  have three sources. More precisely,
\begin{itemize}
\item the inner radicals restrict the domain of analycity inside $D=\cs
\setminus \left(]-\infty, -1] \cup [1, +\infty[\right)$;
\item the pole $z_n=2/(\cos(\sigma/\sqrt n)+\cos(\tau/\sqrt n))$
belongs, for $n$ large enough,
 to $[1, +\infty[$ and does not interfere;
\item the outer radicals do not give further singularities inside,
say, $|z|<2$: indeed, 
\begin{eqnarray*}
2-z(1+e^{-i\sigma/\sqrt n})+2\sqrt{1-z}&=&
2 \cos \theta  e ^{-i\theta} \left( \sqrt{1-z} +1\right)
\left(\sqrt{1-z}+ i \tan \theta \right), \\
%
%
2+z(1-e^{-i\sigma/\sqrt n})+2\sqrt{1+z}&=&
2 \cos \theta  e ^{-i\theta}\left( \sqrt{1+z} +1\right)
\left(1+ i \tan \theta \sqrt{1+z}\right),
%
\end{eqnarray*}
with $\theta = \frac{\sigma}{2\sqrt n}$, and these factorizations
imply that the arguments of
these expressions belong to $]-\pi, \pi[$ as soon as $|z|<2$ and
$\theta < \pi /6$.
\end{itemize}
Consequently, for $\sigma$ and $\tau$ fixed, and $n$ large enough, the
function $S(x,y;z/4)$ is analytic in $D\cap \{|z|<2\}$.
In this domain, the modulus of the numerator of~\Ref{explicite} is bounded by
$10$. The sum of cosines that occurs at the
denominator is larger than $1$ for $n$ large enough, so that finally,
\beq |S(x,y;z/4)| \le \frac{5}{|z_n-z|}.\label{borne-sup}\eeq
 We choose a contour ${\cal C}$ which depends on $n$ and consists of
 four parts ${\cal C}_1$, ${\cal C}_2$, ${\cal C}_3$ and ${\cal
C}_4$ (see Fig.~\ref{hankel}):
\begin{itemize}
\item ${\cal C}_1$ and ${\cal C}_3$ are two symmetric arcs  of radius
$r_n=1+ \log^2 n /n$, centered at the origin,
\item ${\cal C}_2$ a Hankel contour around $1$, at distance
$1/n$ of the real axis, 
which meets  ${\cal C}_1$ and ${\cal C}_3$:
$${\cal C}_2= \left\{ 1+ \frac{t-i}{n}, \ t \in [0, t_n] \right\} \cup
\left\{1-\frac{e^{i\theta}}{n}, \  \theta \in [-\pi/2, \pi/2]\right\}
\cup
\left\{ 1+ \frac{t+i}{n}, \ t \in [0, t_n] \right\} 
$$
where $(1+t_n/n)^2+1/n^2=r_n^2$, so that $t_n \le \log ^2 n$ 
and $t_n=\log^2n+O(1/n)$;
\item ${\cal C}_4=-{\cal C}_2$ is the symmetric Hankel contour around $-1$.
\end{itemize}

\begin{figure}[ht]
\begin{center}
\begin{picture}(0,0)%
\includegraphics{hankel.pstex}%
\end{picture}%
\setlength{\unitlength}{4144sp}%
\begingroup\makeatletter\ifx\SetFigFont\undefined%
\gdef\SetFigFont#1#2#3#4#5{%
  \reset@font\fontsize{#1}{#2pt}%
  \fontfamily{#3}\fontseries{#4}\fontshape{#5}%
  \selectfont}%
\fi\endgroup%
\begin{picture}(4569,3039)(394,-2818)
\put(721,-1591){\makebox(0,0)[lb]{\smash{\SetFigFont{10}{12.0}{\familydefault}{\mddefault}{\updefault}
$-2$
}}}
\put(4546,-1546){\makebox(0,0)[lb]{\smash{\SetFigFont{10}{12.0}{\familydefault}{\mddefault}{\updefault}
$2$
}}}
\put(3646,-1546){\makebox(0,0)[lb]{\smash{\SetFigFont{10}{12.0}{\familydefault}{\mddefault}{\updefault}
$1$
}}}
\put(2386,-1546){\makebox(0,0)[lb]{\smash{\SetFigFont{10}{12.0}{\familydefault}{\mddefault}{\updefault}
$0$
}}}
\put(1531,-1186){\makebox(0,0)[lb]{\smash{\SetFigFont{12}{14.4}{\familydefault}{\mddefault}{\updefault}
${\cal C}_4$
}}}
\put(3466,-1186){\makebox(0,0)[lb]{\smash{\SetFigFont{12}{14.4}{\familydefault}{\mddefault}{\updefault}
${\cal C}_2$
}}}
\put(3556,-2626){\makebox(0,0)[lb]{\smash{\SetFigFont{12}{14.4}{\familydefault}{\mddefault}{\updefault}
${\cal C}_1$
}}}
\put(1441,-466){\makebox(0,0)[lb]{\smash{\SetFigFont{12}{14.4}{\familydefault}{\mddefault}{\updefault}
${\cal C}_3$
}}}
\put(2926,-691){\makebox(0,0)[lb]{\smash{\SetFigFont{12}{14.4}{\familydefault}{\mddefault}{\updefault}
$r_n$
}}}
\end{picture}
\end{center}
\caption{The contour 
$\cal C$.}
\label{hankel}
\end{figure}

The integral~\Ref{cauchy} on $\cal C$ is the sum of the contributions of the
contours ${\cal C}_i$. We shall see that the dominant contribution is
that of ${\cal C}_2$, because of the vicinity of the pole $z_n$.
Let us consider first the contours ${\cal C}_1$
 and ${\cal C}_3$. On these contours, $|z-z_n|>1/n$ so
that by~\Ref{borne-sup},
$
|S(x,y;z/4)|<
5n.
$
Therefore the modulus of the integral on these arcs is 
bounded by 
$
10\pi\,n\,r_n^{-n}
=O(n^{1-\log n})=o(1/n).
$

Consider next the contour ${\cal C}_4$. There, $|z-z_n|>1$, so that
by~\Ref{borne-sup}, $|S(x,y;z/4)|<5$ and the integral is 
small because the contour itself is small: the integral on ${\cal C}_4$
is bounded by $5(2t_n/n+\pi/n)\,(1-1/n)^{-n-1} = O(\log^2n/n)$.
There remains the integral on ${\cal C}_2$.
As $z$ varies along ${\cal C}_2$, the variable $t$ defined by
$z=1+t/n$ varies along the contour 
$\gamma_n$ shown on the left side of Fig.~\ref{gamma-t-infty}.
As $n$ goes to infinity, this contour converges to $\gamma_\infty$. Let
$t\in\gamma_\infty$. Then $t\in\gamma_n$ for $n $ large enough, and,
as $n$ goes to infinity,
the following approximations hold 
with error terms independant of $t$:
\begin{eqnarray*}
\displaystyle \left(2-z(1+\bar x)+2\sqrt{1-z}\right)^{1/2}
&=&\displaystyle {n^{-1/4}}{(2\sqrt{-t}+i\sigma)^{1/2}}
\left(1+O\left(\log^2n/\sqrt n\right)\right)\\
\displaystyle\left(2+z(1-\bar x)+2\sqrt{1+z}\right)^{1/2}
&=&\displaystyle \sqrt2{\textstyle\sqrt{1+\sqrt2}}
\left(1+O(1/\sqrt n)\right)\\
\displaystyle{4-z(x+\bar x+y+\bar y)}
&=&\displaystyle n^{-1}(\sigma^2+\tau^2-4t)
\,\left(1+O\left({\log^2n}/n\right)\right)\\
\displaystyle z^{-n-1} &=&\displaystyle  e^{-t} \left( 1+ O(\log^4 n /n)\right).
\end{eqnarray*}
Observe that in the first approximation, $\Re(2\sqrt{-t}+i\sigma)>0$ for
$t\in\gamma_\infty$,
so that the square root causes no difficulties.
Hence, uniformly in $z \in {\cal C}_2$, we have
$$S(x,y;z/4)z^{-n-1}= n^{3/4}\sqrt2{\textstyle\sqrt{1+\sqrt2}}\,
\frac{(2\sqrt{-t}+i\sigma)^{1/2} e^{-t}} {\sigma^2+\tau^2-4t}
\left( 1+ O(\log^2 n/\sqrt n)\right)$$
with $z=1+t/n$. For $t \in \gamma_n$,  $|\sigma^2+\tau^2-4t|\geq4$ and
$|\exp (-t)|=|\exp(-\Re (t))|\le e$. Moreover, $|t|\le \log ^2 n$, so
that $|i\sigma +2 \sqrt{-t}|=O(\log n)$. Hence the previous identity
implies that
$$S(x,y;z/4)z^{-n-1}= n^{3/4}\sqrt2{\textstyle\sqrt{1+\sqrt2}}\,
\frac{(2\sqrt{-t}+i\sigma)^{1/2} e^{-t}}{\sigma^2+\tau^2-4t}
+ O(n^{1/4} \log^3 n).$$
Let us now integrate this over ${\cal C}_2$:
\begin{eqnarray*}
\displaystyle \int_{{\cal C}_2}S(e^{i\sigma/\sqrt n},e^{i\tau/\sqrt
n};z/4)\frac{dz}{z^{n+1}} &=&
\displaystyle\frac{1}{n}\ \int_{{\gamma}_n}S(e^{i\sigma/\sqrt n},e^{i\tau/\sqrt
n};(1+t/n)/4)\frac{dt}{(1+t/n)^{n+1}}\\
&=& \displaystyle n^{-1/4}\sqrt2{\textstyle\sqrt{1+\sqrt2}}\,
\int_{{\gamma}_n}\frac{(2\sqrt{-t}+i\sigma)^{1/2}e^{-t}}{\sigma^2+\tau^2-4t}
dt
+O(n^{-3/4} \log^5 n).
\end{eqnarray*}
 As $n$ goes to infinity,
$$  \int_{\gamma_n}
\frac{(2\sqrt{-t}+i\sigma)^{1/2}e^{-t}}{\sigma^2+\tau^2-4t} dt
 \longrightarrow \int_{\gamma_\infty}
\frac{(2\sqrt{-t}+i\sigma)^{1/2}e^{-t}}{\sigma^2+\tau^2-4t} dt .$$ 
Hence finally,
$$ \int_{{\cal C}_2}S(e^{i\sigma/\sqrt n},e^{i\tau/\sqrt
n};z/4)\frac{dz}{z^{n+1}} = n^{-1/4}\sqrt2{\textstyle\sqrt{1+\sqrt2}}\,
\int_{{\gamma}_\infty}\frac{(2\sqrt{-t}+i\sigma)^{1/2}e^{-t}}{\sigma^2+\tau^2-4t}dt
\ \left(1 + o(1)\right).$$
Thus ${\cal C}_2$ is really the part of the contour that yields the
dominant contribution to the integral of Eq.~\Ref{cauchy}. We now
inject in~\Ref{cauchy} the following ingredients:

--  the four estimates of the
integrals on the contours ${\cal C}_i$,

--  the asymptotic behaviour of
$a(n)$, which follows Theorem~\ref{dim2},

--   the complement formula, according to which $\Gamma(1/4)
    \Gamma(3/4)= \sqrt 2 \pi$.

\noindent Proposition~\ref{charact} follows.
\cqfd

\noindent\textbf{Proof of Theorem \ref{distribXY}.}  There remains to
check that the limit function $\Phi$, given in
%
%
Proposition~\ref{charact}, coincides with the
characteristic function $\Psi$ of the  distribution
defined in Theorem~\ref{distribXY}. We use the expression of the density
in polar coordinates to express this 
characteristic function:
$$\Psi(\sigma,\tau) = \GE\left( e^{i\sigma X +i\tau Y}\right)
=\frac{2}{\Gamma(1/4)}
\int_{-\pi}^\pi d\theta
\int_0^\infty d\rho\;
e^{i\rho(\sigma \cos\theta+\tau \sin\theta)}\;
\rho^{3/2}e^{-\rho^2}\cos(\theta/2).$$
One possible approach is to expand $\Psi(\sigma,\tau)$ and
$\Phi(\sigma,\tau)$ in series of $\sigma$ and $\tau$, and to check that
the coefficients coincide (these coefficients are, essentially, the
moments of the pair $(X,Y)$). This  natural
approach works, but requires a few
more calculations 
than the method we present below.

We choose to work with polar coordinates, not only for the density of
$(X,Y)$, but also for the variables $\sigma$ and $\tau$, which we 
take to be $r\cos\phi$ and $r\sin\phi$ respectively.
 The characteristic function $\Psi$ becomes
$$
\Psi(\sigma,\tau) =\Psi(r\cos\phi,r\sin\phi)=
\frac{2}{\Gamma(1/4)}
\int_{-\pi}^\pi d\theta
\int_0^\infty d\rho\;
e^{ir\rho\cos(\theta-\phi)}\;
\rho^{3/2}e^{-\rho^2}\cos(\theta/2).
$$
Let us expand the integrand in $r$:
\[
\Psi(\sigma,\tau) =
\frac{2}{\Gamma(1/4)}
\int_{-\pi}^\pi d\theta
\int_0^\infty d\rho
\sum_{n\geq0}\frac{(ir)^n}{n!}
\rho^n\rho^{3/2}e^{-\rho^2}
 \cos(\theta-\phi)^n\cos(\theta/2).
\]
This (triple) sum is absolutely convergent so that we can
exchange the sum and the integrals, and then separate the integrals on
$\rho$ and $\theta$. Using the
definition of the Gamma function, 
\[
\Gamma(s)=\int_0^{\infty}e^{-t}t^{s-1} dt
=2 \int_0^{\infty}\rho^{2s-1}e^{-\rho^2}d\rho,
\]
we then  evaluate the integral over $\rho$, and  obtain:
\[
\Psi(r\cos\phi,r\sin\phi)=
\frac{1}{\Gamma(1/4)}
\sum_{n\geq0}\frac{(ir)^n}{n!}
\Gamma\left(\frac n2+\frac54\right)
\int_{-\pi}^\pi d\theta
\cos(\theta-\phi)^n \cos(\theta/2).
\]
Our aim is to transform this expression into the expression of
Proposition~\ref{charact}. In particular, we need 
to introduce the Hankel contour $\gamma_\infty$, which is 
known to occur in the Hankel representation of the {\em inverse\/} of
the Gamma function:
\beq
\frac1{\Gamma(s)}=\frac1{2\pi i}\int_{\gamma_\infty}(-t)^{-s}e^{-t}dt.
\label{gamma-inverse}
\eeq
The trouble is that the Gamma function appears only as a numerator,
and not as a denominator, in the expression of $\Psi$. We shall remedy this
thanks to 
the duplication formula:
\[2^{2s-1} \Gamma(s) \Gamma(s+1/2) = \sqrt \pi \, \Gamma(2s).
\]
Applied to $s=n/2 + 3/4$, it
allows us to rewrite $\Psi(\sigma,\tau)$ as
\[
\Psi(\sigma,\tau)=
\frac{\sqrt{\pi}}{\sqrt 2\Gamma(1/4)}
\sum_{n\geq0}\frac{(ir)^n}{2^n n!}
\frac{\Gamma\left(n+\frac32\right)}
{\Gamma\left(\frac n2+\frac34\right)}
\int_{-\pi}^\pi d\theta
\cos(\theta-\phi)^n \cos(\theta/2)
\]
and to introduce, at last, the contour $\gamma_\infty$:
\[
\Psi(\sigma,\tau)=
\frac{1}{2i \sqrt{2\pi}\Gamma(1/4)}
\sum_{n\geq0}\left(\frac{ir}{2}\right)^n
\frac{\Gamma\left(n+\frac32\right)}
{n!}
\int_{\gamma_\infty}dt(-t)^{-\frac n2-\frac34}\,e^{-t}
\int_{-\pi}^\pi d\theta
\cos(\theta-\phi)^n\cos(\theta/2) .
\]
But the summation over $n$ is now 
subject to an explicit resummation. Indeed, for $|z|<1$,
\beq
\sum_{n\geq0}z^n\frac{\Gamma\left(n+\frac32\right)}{n!}
=
\frac{\sqrt{\pi}}{2(1-z)^{3/2}}. \label{3/2} \eeq
Upon exchanging the sum and the integrals in the above expression of
$\Psi(\sigma, \tau)$   we obtain:
\[
\Psi(\sigma,\tau)=
\frac{1}{2i \sqrt{2\pi}\Gamma(1/4)}
\int_{\gamma_\infty}dt(-t)^{-\frac34}\,e^{-t}
\int_{-\pi}^\pi d\theta
\cos(\theta/2)
\sum_{n\geq0}\left(\frac{ir\cos(\theta-\phi)}{2\sqrt{-t}}\right)^n
\frac{\Gamma\left(n+\frac32\right)}
{n!}.
\]
As $|t|\ge 1$,  this exchange of summations is valid if $r< 2$ by
virtue of the 
absolute convergence of the power series~\Ref{3/2}.
 For larger values of $r$, it suffices to replace the Hankel contour
$\gamma_\infty$ by $2r^2\gamma_\infty$ in~\Ref{gamma-inverse} and the
above lines to obtain an absolutely convergent series. Finally, 
\beq
\Psi(\sigma,\tau )=
\frac{1}{4i\sqrt 2 \Gamma(1/4)}
\int_{\gamma_\infty}dt(-t)^{-\frac34}\,e^{-t}
\int_{-\pi}^\pi d\theta
\frac{\cos(\theta/2)}
{(1-ir({-t})^{-1/2}\cos(\theta-\phi)/2)^{3/2}}.
\label{psi-presque}\eeq
 The integration on $\theta$ can be  performed easily, because the integrand has an explicit primitive:
\[
\frac{\cos(\theta/2)}{(1-z\cos(\theta-\phi))^{3/2}}=
\frac2{1-z^2}\frac{\partial}{\partial\theta}
 \left( \frac{\sin(\theta/2) +z \sin(\theta /2 -\phi)}
{(1-z\cos(\theta-\phi))^{1/2}}\right).
\]
Using this primitive, we obtain
\[
\int_{-\pi}^\pi d\theta
\frac{\cos(\theta/2)}{(1-z\cos(\theta-\phi))^{3/2}}
=
\frac{4(1+z\cos\phi)^{1/2}}{1-z^2}.
\]
Applying the case  $z=ir(-t)^{-1/2}/2$ of this identity to
Eq.~\Ref{psi-presque} gives, after a few reductions,
\[
\Psi(\sigma,\tau)=
\frac{2}{i\Gamma(1/4)}
\int_{\gamma_\infty}dt\,e^{-t}
\frac{(2\sqrt{-t}+ir\cos\phi)^{1/2}}
{r^2-4t}= \Phi(r\cos\phi,r\sin\phi).
\]
\cqfd

\noindent
{\bf Acknowledgements.} First of all, we are greatly indebted to
Olivier Roques, who 
discovered the remarkable conjectures that were the starting point of
this work, and to Philippe ``Duduche'' Duchon, who brought
them to our attention.  
Then, we have to say that the current
version of this paper has little to do with the very first version,
which we wrote almost two years ago. The differences
stem party from simplifications in the proof of the main result, but
mostly from the numerous discussions or mail exchanges we had with
several colleagues; some of them had a definite influence on the form
and content of this paper, like  B\'etr\'ema,  Bertrand
Duplantier, Ira Gessel, 
Barry Hughes, 
and Jean-Fran\c cois Le Gall.
MBM also had many interesting discussions with the attendants of the
 workshop ``Self-interacting Random Processes''  hold in Oberwolfach,
 in May 2000. Finally, we  thank for their
 interest and patience our 
colleagues from Bordeaux, Melbourne and Nancy, which we have kept
 bothering with this topic  for almost
two years.


\begin{thebibliography}{99}


\bibitem{abhyankar} S. S. Abhyankar, {\em Algebraic geometry for scientists
 and engineers}, Mathematical surveys and monographs {\bf 35}, American
 Mathematical Society, 1990.


\bibitem{6gus} C. Banderier,  M. Bousquet-M\'elou,
 A. Denise, P. Flajolet, D. Gardy and D.
 Gouyou-Beauchamps,  Generating functions for generating trees,
 to appear in {\em Discrete Math.}

\bibitem{florence} E. Barcucci, E. Pergola, R. Pinzani and S. Rinaldi,
A bijection for some paths on the slit plane, preprint 2000,
%
%
Universit\`a di Firenze.


\bibitem{prep} M. Bousquet-M\'elou, Walks on the slit plane: other
approaches, in preparation. 

\bibitem{BoPe98}
M.~Bousquet-M\'elou and M.~Petkov\v{s}ek,
\newblock Linear recurrences with constant coefficients: the multivariate case,
{\em Discrete Math.} {\bf 225} (2000) 51--75.

\bibitem{considine}
D. Considine and S. Redner, Repulsion of random and self-avoiding
walks from excluded points and lines, {\em J. Phys. A: Math. Gen.} {\bf
22} (1989) 1621--1638.

\bibitem{pizaler} M.-P. Delest and G. Viennot, 
Algebraic languages and polyominoes enumeration, 
{\em Theoret. Comput. Sci.} {\bf 34} (1984) 169--206.

\bibitem{FaIa}
G.~Fayolle and R.~Iasnogorodski,
\newblock Solutions of functional equations arising in the analysis of
  two-server queueing models.
\newblock In {\em Performance of computer systems}, pages 289--303.
  North-Holland, 1979.

\bibitem{Malyshev} G. Fayolle, R. Iasnogorodski and V.
Malyshev, {\em Random walks in the quarter plane : algebraic methods, boundary
value problems, and applications}.  
Applications of Mathematics {\bf 40}, Springer, New York, 1999.


\bibitem{Flajolet87}
P.~Flajolet,
\newblock Analytic models and ambiguity of context--free languages,
\newblock {\em Theoret. Comput. Sci.} {\bf 49} (1987) 283--309.


\bibitem{flajolet-odlyzko}
P. Flajolet and A. Odlyzko, Singularity analysis of generating
functions, {\em SIAM J. Disc. Math.} {\bf 3} No. 2 (1990) 216--240.

\bibitem{gessel} I. Gessel, personnal communication, September 2000.

\bibitem{gj} I. P. Goulden and D. M. Jackson, {\em Combinatorial
 enumeration,} John Wiley and  Sons, 1983.


\bibitem{kesten} H. Kesten, Hitting probabilities of random walks on
$\zs ^d$, {\em Stoch. Proc. and Appl.} {\bf 25} (1987) 165--184.


\bibitem{knuth} D.~E.~Knuth,
   {\em The Art of Computer Programming, Vol.~$1$: Fundamental Algorithms}.
   Addison-Wesley, Reading Mass., 1968.
\newblock Exercises 4 and 11, Section 2.2.1.


\bibitem{lawler} G. F. Lawler, {\em Intersections of random walks},
Probabilities and its applications, Birkh\"auser Boston, 1991.

\bibitem{legall} J.-F. Le Gall,  personnal communication, May 2000.




\bibitem{ray} D. Ray, Stable processes with an absorbing barrier, {\em
Trans. Amer. Math. Soc.} {\bf 89 } (1958) 16--24.


\bibitem{sloane} N. J. A. Sloane and S. Plouffe, {\em The encyclopedia of
integer sequences}, Academic Press, 1995.


\end{thebibliography}
\end{document}